\documentclass[10pt, a4paper]{article}
\usepackage{cmap}
\usepackage[main=english,russian]{babel}
\usepackage[utf8]{inputenc}
\usepackage[text={8in,10in},centering]{geometry}
\usepackage{amsfonts,amsmath,amsxtra,amsthm,amssymb,latexsym, tabularray}
\usepackage{graphicx}
\usepackage{verbatim}
\usepackage{cite, color}
\usepackage{subcaption}
\usepackage{float}
\usepackage{comment}
\theoremstyle{plain}
\newtheorem{theorem}{Theorem}[section]

\newtheorem{lemma}{\bf \indent Lemma}[section]
\newtheorem{corollary}{\bf \indent Corollary}[section]
\newtheorem*{theorem*}{\bf \indent Theorem}

\theoremstyle{remark}

\newtheorem{definition}{\bf \indent Definition}[section]
\newtheorem{example}{\bf \indent Example}[section]

\numberwithin{equation}{section}

\def\R{{\mathbb R}}

\newcommand{\VEC}{\operatorname{Vec}\nolimits}

\newcommand{\be}[1]{\begin{equation}\label{#1}}
\newcommand{\ee}{\end{equation}}

\begin{document}
\parindent=0cm

\title{Sectional Curvature and Structural Functions of a Two-Dimensional or Three-Dimensional Lorentzian Manifold
}
\author{A.Z. Ali,\ Yu.L. Sachkov}
\date{}

\maketitle

\section{Introduction}
In this note, we derive a formula for the sectional curvature on an arbitrary two-dimensional or three-dimensional smooth manifold $M$ equipped with a Lorentzian metric $g$ (a symmetric non-degenerate smooth tensor field of second order of index one). Here we assume that $g$ satisfies the following property: there exist basis vector fields $X_1$, $X_2$ on $M$ such that $g(X_1,X_1) = a_{11}$, $g(X_2,X_2) = a_{22}$, $g(X_1,X_2) = a_{12} = g(X_2,X_1) = a_{21}$, where $a_{ij}$, $i,j=1,2$ are constant real numbers.

The note provides a scheme for computing the sectional curvature, the derivation of a general formula, as well as some examples of the application of the obtained result.

\section{Problem Statement}
Let us recall the necessary definitions and known statements \cite{beem}, \cite{lor_lob}.

Let $M$ be a smooth manifold. Denote by $Vec(M)$ the set of smooth vector fields on $M$.

\begin{definition}
A connection $D$ on a manifold $M$ is a mapping $D: (\VEC(M))^2 \rightarrow \VEC(M)$ satisfying the axioms:
\begin{itemize}
    \item[$(1)$] $D_{fU+gV}W = fD_{U}W+gD_{V}W$ for any $U$, $V$, $W$ $\in$ $\VEC(M)$, for any $f, g \in C^{\infty}(M)$,
    \item[$(2)$] $D_{V}(\alpha_1 W_1 + \alpha_2 W_2) = \alpha_1 D_{V}W_1 + \alpha_2 D_{V}W_2$ for any $W_1$, $W_2$ $\in \VEC(M)$, for any $\alpha_1$, $\alpha_2$ $\in \mathbb{R}$,
    \item[$(3)$] $D_{V}(fW) = (Vf)W + fD_{V}W$ for any $V$, $W$ $\in \VEC(M)$, for any $f \in C^{\infty}(M)$.
\end{itemize}
\end{definition}

The vector field $D_{V}W$ is called the {\it covariant derivative of the field $W$ along $V$ for the connection $D$}.

\begin{theorem}
On a Lorentzian manifold $(M, g)$, there exists a unique connection $D$ such that:
\begin{itemize}
    \item[$(4)$] $[V, W] = D_{V}W - D_{W}V$,
    \item[$(5)$] $Xg(V, W) = g(D_X V, W) + g(V, D_{X} W)$,
\end{itemize}
for any $X, V, W \in \VEC(M)$. $D$ is called the {\it Levi-Civita connection on $M$}, and is characterized by the Koszul formula
\begin{equation}
\label{koszul}
2g(D_{V}W,X) =Vg(W,X)+Wg(X,V)-Xg(V,W)-g(V,[W,X])+g(W,[X,V])+g(X,[V,W]).
\end{equation}
\end{theorem}

The mapping $R: (\VEC(M))^3 \rightarrow C^{\infty}(M)$, given by the formula $R_{XY}Z = D_{[X,Y]}Z-[D_{X},D_{Y}]Z$, where $X, Y, Z \in \VEC(M)$, is called the {\it Riemann curvature tensor of $(M,g)$}.

Let $q \in M$, $P$ be a two-dimensional plane in the tangent space $T_{q}M$. Then for vectors $v, w \in T_{q}M$, define $Q(v,w) = g(v,v)g(w,w)-(g(v,w))^2$. The plane $P$ is called {\it non-degenerate} if $Q(v,w) \neq 0$ for some (and then for any) basis $v, w$ of the plane $P$.

\begin{lemma}
Let $P \subseteq T_{q}M$ be a non-degenerate two-dimensional plane. Then the number
\begin{equation}
\label{sectionalcurvature}
    K(q,P) = \frac{g(R_{vw}v,w)}{Q(v,w)}
\end{equation}
does not depend on the choice of basis in the plane $P$, and is called the {\bf sectional curvature of the plane sector $P$}.
\end{lemma}

\begin{theorem}
Let $(M,g)$ be a Lorentzian manifold of dimension $n \geqslant 2$, and let $K \in \mathbb{R}$. Then the following conditions are equivalent:
\begin{itemize}
    \item[$(1)$] $M$ has constant curvature $K$,
    \item[$(2)$] for any point $q \in M$, there exists a neighborhood isometric to an open subset of de Sitter space $\mathbb{S}_1^n$ for $K>0$, Minkowski space $\mathbb{R}_{1}^{n}$ for $K = 0$, or anti-de Sitter space $\widetilde{\mathbb{H}_{1}^{n}}$ for $K < 0$.
\end{itemize}
\end{theorem}

\section{Two-Dimensional Manifold}

\subsection{Plan for Computing the Sectional Curvature for a Two-Dimensional Manifold}

\begin{itemize}
    \item[$(1)$] Compute the commutator $[X_1, X_2]$;
    \item[$(2)$] Compute the Levi-Civita connection coefficients using the Koszul formula (\ref{koszul}), using the basis vector fields $X_1, X_2$;
    \item[$(3)$] Compute the Riemann curvature tensor $R_{v w}v$, where $v = X_1$, $w = X_2$;
    \item[$(4)$] Compute $Q(v,w)$, where $v = X_1$, $w = X_2$;
    \item[$(5)$] Obtain the sectional curvature $K$ using the performed computations and formula (\ref{sectionalcurvature}).
\end{itemize}

\subsection{Main Formula, Useful Corollaries, and Examples}

\begin{theorem}
The sectional curvature of a Lorentzian two-dimensional manifold $M$ with Lorentzian metric $g$, its basis vector fields $X_1$, $X_2$, $[X_1,X_2] = c_{12}^{1}X_1+c_{12}^{2}X_2$, and constant metric values on them: $g(X_1,X_1) = a_{11}$, $g(X_2,X_2) = a_{22}$, $g(X_1,X_2) = a_{12} = g(X_2,X_1) = a_{21}$, equals
\begin{equation}
\label{dvum_kriv_obsh}
    K= \frac{a_{12}\xi_1 + a_{22}\xi_2}{\det{\mathcal{A}}^2},
\end{equation}
where
\begin{equation*}
    \xi_1 = a_{12} \left[ \left( c_{12}^1A + c_{12}^2B \right) - \left( X_1 B \right) + \left( X_2 A \right) \right], \quad \xi_2 = a_{11} \left[ - \left( c_{12}^1A + c_{12}^2B \right) + \left( X_1B \right) - \left( X_2A \right) \right],
\end{equation*}
\begin{equation*}
A = a_{11}c_{12}^1+a_{12}c_{12}^2, \quad B = a_{12}c_{12}^1+a_{22}c_{12}^2,
\end{equation*}

\begin{equation}
\label{dvum_matr}
    \mathcal{A} = 
    \begin{pmatrix}
    a_{11} & a_{12}\\
    a_{12} & a_{22}
    \end{pmatrix}. 
\end{equation}

\end{theorem}

\begin{corollary}
For orthogonal $X_1$, $X_2$, i.e. $a_{12} = a_{21} = 0$, $a_{11} \neq 0$, $a_{22} \neq 0$, the sectional curvature equals
\begin{equation}
\label{dvum_kriv_ortogon}
K = \frac{-a_{11}\left(c_{12}^1\right)^2 -a_{22} \left(c_{12}^2\right)^2 - a_{11} \left(X_2 c_{12}^1\right) + a_{22} \left(X_1 c_{12}^2\right)}{a_{11}}.
\end{equation}
\end{corollary}

\begin{corollary}
For orthonormal $X_1$, $X_2$, i.e. $a_{11} = -1$, $a_{12} = a_{21} = 0$, $a_{22} = 1$, the sectional curvature equals
\begin{equation}
\label{dvum_kriv_ortonorm}
K = -\left(c_{12}^1\right)^2 + \left(c_{12}^2\right)^2 - \left(X_2 c_{12}^1\right) - \left(X_1 c_{12}^2\right).
\end{equation}
\end{corollary}

\begin{example}
{\bf Two-dimensional anti-de Sitter space $\widetilde{\mathbb{H}_{1}^{2}}$.}

The Lorentzian metric has the expression: $g =  - \ch^2{\theta} d \varphi^2+d \theta^2$, and the orthonormal frame: $X_1 = \frac{1}{\ch{\theta}} \frac{\partial}{\partial \varphi},\ X_2 = \frac{\partial }{\partial \theta}$. We have $g(X_1,X_1) = -1$, $g(X_2,X_2) = 1$, $g(X_1,X_2) = g(X_2,X_1) = 0$, and also $[X_1,X_2] = \th{\theta}X_1$, i.e. the structural functions are as follows: $c_{12}^1 = \th{\theta}$, $c_{12}^2 = 0$.

Therefore, by formula (\ref{dvum_kriv_ortonorm}), the sectional curvature of two-dimensional anti-de Sitter space equals:
\begin{equation*}
K = -\left[\left(c_{12}^1\right)^2 - \left(c_{12}^2\right)^2 + \left(X_2 c_{12}^1\right) + \left(X_1 c_{12}^2\right)\right] = -\left[ (\th{\theta})^2 + \partial_{\theta}(\th{\theta}) \right] = -\frac{\sh^2{\theta}+1}{\ch^2{\theta}} = -1,
\end{equation*}
which coincides with the result obtained in \cite{lor_lob}.
\end{example}

\begin{example}
\label{igreki}
Let $X_1$, $X_2$ be a basis of a two-dimensional solvable Lie algebra such that $[X_1,X_2] = -X_1$, and $Y_1$, $Y_2$ be an orthonormal basis of the Lorentzian metric $g$, with $Y_1 = \alpha X_1 + \gamma X_2$, $Y_2 = \beta X_1 + \delta X_2$.

Let us compute the sectional curvature using formula (\ref{dvum_kriv_ortonorm}). To do this, we first need to compute the structural functions.
\begin{equation*}
[Y_1,Y_2] = [\alpha X_1 + \gamma X_2,\beta X_1 + \delta X_2] = \alpha \beta[X_1,X_1] + \alpha \delta [X_1,X_2] + \gamma \beta [X_2,X_1] + \gamma \delta [X_2,X_2] = (\gamma \beta - \alpha \delta)X_1.
\end{equation*}
Expand the obtained vector in terms of $Y_1$, $Y_2$:
\begin{equation*}
\gamma Y_2 - \delta Y_1 = \gamma( \beta X_1 + \delta X_2 ) - \delta ( \alpha X_1 + \gamma X_2 ) = (\gamma \beta - \alpha \delta)X_1,
\end{equation*}
therefore
\begin{equation*}
[Y_1,Y_2] = \gamma Y_2 - \delta Y_1 = c_{12}^1 Y_1 + c_{12}^2 Y_2,
\end{equation*}
i.e. $c_{12}^1 = -\delta$, $c_{12}^2 = \gamma$.
We obtain:
\begin{equation*}
K = - \left[ \left(c_{12}^1\right)^2 - \left(c_{12}^2\right)^2 + \left(Y_2 c_{12}^1\right) + \left(Y_1 c_{12}^2\right) \right] = - \left[ \delta^2 - \gamma^2 \right] = \gamma^2 - \delta^2,
\end{equation*}
which coincides with the result obtained in \cite{lor_lob}.
\end{example}

\begin{example}
{\bf Left-invariant problems on a two-dimensional Lie group $G$.}

Let $X_1$, $X_2$ be vector fields such that $[X_2,X_1] = X_1$, i.e. $c_{12}^1 = -1$, $c_{12}^2 = 0$ in our notation, and $g(X_1,X_1) = c^2-a^2 = a_{11}$, $g(X_1,X_2) = g(X_2, X_1) = cd-ab = a_{12}$, $g(X_2,X_2) = d^2-b^2 = a_{22}$.

Let us compute the sectional curvature using formula (\ref{dvum_kriv_ortogon}).

Since the structural functions are constants, all Lie derivatives along our vector fields are zero. Therefore, it suffices to compute $A$ and $B$ to find $\xi_1$, $\xi_2$.
\begin{equation*}
A = a_{11}c_{12}^1 + a_{12}c_{12}^2 = -a_{11} = a^2-c^2,\ B = a_{12}c_{12}^1+a_{22}c_{12}^2 = -a_{12}= ab-cd,
\end{equation*}
therefore
\begin{equation*}
\xi_1 = a_{12} \left[ \left( c_{12}^1A + c_{12}^2B \right) - \left( X_1 B \right) + \left( X_2 A \right) \right] = -a_{12}A = -a_{12}( a_{11}c_{12}^1 + a_{12}c_{12}^2 ) = a_{11}a_{12},
\end{equation*}
\begin{equation*}
\xi_2 =  a_{11} \left[ - \left( c_{12}^1A + c_{12}^2B \right) + \left( X_1B \right) - \left( X_2A \right) \right] = a_{11}( a_{11}c_{12}^1 + a_{12}c_{12}^2) = -(a_{11})^2.
\end{equation*}
\begin{equation*}
a_{12}\xi_1 + a_{22}\xi_2 = a_{12}a_{11}a_{12} - a_{22}(a_{11})^2 = -a_{11}\left( a_{11}a_{22} - (a_{12})^2 \right) = -a_{11}\det{\mathcal{A}}.
\end{equation*}
We also note that
\begin{equation*}
\det{\mathcal{A}} = a_{11}a_{22}-(a_{12})^2 = (c^2-a^2)(d^2-b^2)-(cd-ab)^2 = -(ad-bc)^2 = -\det(\mathcal{B})^2,
\end{equation*}
where
\begin{equation*}
\mathcal{B} = 
\begin{pmatrix}
a & b \\
c & d
\end{pmatrix}.
\end{equation*}
Then
\begin{equation*}
K = \frac{a_{12}\xi_1 + a_{22}\xi_2}{\det{\mathcal{A}}^2} = \frac{-a_{11}\det{\mathcal{A}}}{\det(\mathcal{B})^4} = \frac{a_{11}}{\det(\mathcal{B})^2},
\end{equation*}

which coincides with the result obtained in \cite{lor_lob}.

\end{example}

\subsection{Derivation of the Main Formula}

We follow the plan.

\begin{itemize}
    \item[$(1)$] We have written the commutator in general form:
    \begin{equation*}
        [X_1,X_2] = c_{12}^{1}X_1+c_{12}^{2}X_2.
    \end{equation*}
    \item[$(2)$] Since for $X_1, X_2$, $a_{ij} = const$, $\forall$ $i,j = 1,2$, the first three terms on the right-hand side of the Koszul formula (\ref{koszul}) are zero if $X_1, X_2$ are substituted for $V, W, X$ in any combinations. Therefore, only the fourth, fifth, and sixth elements of the right-hand side of formula (\ref{koszul}) are meaningful. In total, we need to compute 8 coefficients, namely $g(D_{X_1}X_1,X_i),\ i=1,2$; $g(D_{X_1}X_2,X_i),\ i=1,2$; $g(D_{X_2}X_1,X_i),\ i=1,2$; $g(D_{X_2}X_2,X_i),\ i=1,2$. Let's proceed with the calculations. We can immediately say that $g(D_{X_1}X_1,X_1) = 0$, $g(D_{X_2}X_2,X_2) = 0$, since all commutators consist of $[X_1,X_1]$ or $[X_2,X_2]$ respectively, which are identically zero.
    \begin{equation*}
    g(D_{X_1}X_1,X_2) = \frac{1}{2}(-g(X_1,[X_1,X_2])+g(X_1,[X_2,X_1])+g(X_2,[X_1,X_1]))= 
    \end{equation*}
    \begin{equation*}
    = \frac{1}{2}( -g(X_1,c_{12}^{1}X_1+c_{12}^{2}X_2) + g(X_1,-c_{12}^{1}X_1-c_{12}^{2}X_2)) =
    \end{equation*}
    \begin{equation*}
    = \frac{1}{2}( -a_{11}c_{12}^1-a_{12}c_{12}^2-a_{11}c_{12}^1-a_{12}c_{12}^2 ) = -a_{11}c_{12}^1-a_{12}c_{12}^2 = -A,
    \end{equation*}
    \begin{equation*}
    g(D_{X_1}X_2,X_1) = \frac{1}{2}(-g(X_1,[X_2,X_1])+g(X_2,[X_1,X_1])+g(X_1,[X_1,X_2])) = 
    \end{equation*}
    \begin{equation*}
    = \frac{1}{2}(-g(X_1,-c_{12}^{1}X_1-c_{12}^{2}X_2)+g(X_1,c_{12}^{1}X_1+c_{12}^{2}X_2)) =
    \end{equation*}
    \begin{equation*}
    = \frac{1}{2}( a_{11}c_{12}^1 + a_{12}c_{12}^2 + a_{11}c_{12}^1 + a_{12}c_{12}^2 ) = a_{11}c_{12}^1 + a_{12}c_{12}^2 = A,
    \end{equation*}
    \begin{equation*}
    g(D_{X_1}X_2,X_2) = \frac{1}{2}( -g(X_1,[X_2,X_2])+g(X_2,[X_2,X_1])+g(X_2,[X_1,X_2]) ) = \end{equation*}
    \begin{equation*}
    = \frac{1}{2}( g(X_2,-c_{12}^{1}X_1-c_{12}^{2}X_2)+g(X_2,c_{12}^{1}X_1+c_{12}^{2}X_2) ) = 0,
    \end{equation*}
    \begin{equation*}
    g(D_{X_2}X_1,X_1) = \frac{1}{2}( -g(X_2,[X_1,X_1])+g(X_1,[X_1,X_2])+g(X_1,[X_2,X_1]) ) =
    \end{equation*}
    \begin{equation*}
    =  \frac{1}{2}( g(X_1,c_{12}^{1}X_1+c_{12}^{2}X_2) + g(X_1,-c_{12}^{1}X_1-c_{12}^{2}X_2) ) = 0,
    \end{equation*}
    \begin{equation*} g(D_{X_2}X_1,X_2) = \frac{1}{2}( -g(X_2,[X_1,X_2])+g(X_1,[X_2,X_2])+g(X_2,[X_2,X_1]) ) =
    \end{equation*}
    \begin{equation*}
    = \frac{1}{2}( -g(X_2,c_{12}^{1}X_1+c_{12}^{2}X_2) + g(X_2,-c_{12}^{1}X_1-c_{12}^{2}X_2) ) =
    \end{equation*}
    \begin{equation*} 
    = \frac{1}{2}( -a_{12}c_{12}^1-a_{22}c_{12}^2-a_{12}c_{12}^1 - a_{22}c_{12}^2 ) = -a_{12}c_{12}^1-a_{22}c_{12}^2 = -B,
    \end{equation*}
    \begin{equation*} g(D_{X_2}X_2,X_1) = \frac{1}{2}( -g(X_2,[X_2,X_1])+g(X_2,[X_1,X_2])+g(X_1,[X_2,X_2]) ) = 
    \end{equation*}
    \begin{equation*}
    = \frac{1}{2}( -g(X_2,-c_{12}^{1}X_1-c_{12}^{2}X_2)+g(X_2,c_{12}^{1}X_1+c_{12}^{2}X_2) ) 
    \end{equation*}
    \begin{equation*} 
    = \frac{1}{2}( a_{12}c_{12}^1+a_{22}c_{12}^2 + a_{12}c_{12}^1 + a_{22}c_{12}^2 ) = a_{12}c_{12}^1+a_{22}c_{12}^2 = B.
    \end{equation*}

    Now let's compute the covariant derivatives.

    \begin{itemize}
        \item Temporarily denote
        \begin{equation*}
            D_{X_1}X_1 = uX_1 + vX_2.
        \end{equation*}
        Find $u$ and $v$ from the system of equations:
        \begin{equation*}
        \begin{cases}
            g\left( D_{X_1}X_1,X_1 \right) = 0,\\
            g\left(D_{X_1}X_1,X_2 \right) = -A.
        \end{cases}
        \end{equation*}
        \begin{equation*}
        \begin{cases}
            g\left( uX_1 + vX_2,X_1 \right) = 0,\\
            g\left(uX_1 + vX_2,X_2 \right) = -A,
        \end{cases}
        \Leftrightarrow
        \begin{cases}
    ug\left(X_1,X_1\right)+vg\left(X_2,X_1 \right) = 0,\\
    ug\left(X_1,X_2\right)+vg\left(X_2,X_2 \right) = -A,
        \end{cases}
        \Leftrightarrow
        \begin{cases}
            ua_{11} + va_{12} = 0,\\
            ua_{12} + va_{22} = -A,
        \end{cases}
        \Leftrightarrow
        \end{equation*}
    \begin{equation*}
    \Leftrightarrow
    \begin{cases}
        u = \frac{a_{12}A}{a_{11}a_{22}-(a_{12})^2} = \frac{a_{12}A}{\det{(\mathcal{A})}},\\
        v = \frac{-Aa_{11}}{a_{11}a_{22}-(a_{12})^2} = \frac{-Aa_{11}}{\det{(\mathcal{A})}}.
    \end{cases}
    \end{equation*}

    \item Temporarily denote
        \begin{equation*}
            D_{X_1}X_2 = uX_1 + vX_2.
        \end{equation*}
        Find $u$ and $v$ from the system of equations:
        \begin{equation*}
        \begin{cases}
            g\left( D_{X_1}X_2,X_1 \right) = A,\\
            g\left(D_{X_1}X_2,X_2 \right) = 0.
        \end{cases}
        \end{equation*}
        \begin{equation*}
        \begin{cases}
        g\left( uX_1 + vX_2,X_1 \right) = A,\\
        g\left(uX_1 + vX_2,X_2 \right) = 0.
        \end{cases}
        \Leftrightarrow
        \begin{cases}
        ua_{11} + va_{12} = A,\\
        ua_{12} + va_{22} = 0,
        \end{cases}
        \Leftrightarrow
        \begin{cases}
        u = \frac{a_{22}A}{\det{(\mathcal{A})}},\\
        v = \frac{-a_{12}A}{\det{(\mathcal{A})}}.
        \end{cases}
        \end{equation*}

        \item Temporarily denote
        \begin{equation*}
            D_{X_2}X_1 = uX_1 + vX_2.
        \end{equation*}
        Find $u$ and $v$ from the system of equations:
        \begin{equation*}
        \begin{cases}
            g\left( D_{X_2}X_1,X_1 \right) = 0,\\
            g\left(D_{X_2}X_1,X_2 \right) = -B.
        \end{cases}
        \end{equation*}
        \begin{equation*}
        \begin{cases}
        g\left( uX_1 + vX_2,X_1 \right) = 0,\\
            g\left(uX_1 + vX_2,X_2 \right) = -B,
        \end{cases}
        \Leftrightarrow
        \begin{cases}
        ua_{11} + va_{12} = 0,\\
        ua_{12} + va_{22} = -B,
        \end{cases}
        \Leftrightarrow
        \begin{cases}
        u = \frac{a_{12}B}{\det{(\mathcal{A})}},\\
        v = \frac{-a_{11}B}{\det{(\mathcal{A})}}.
        \end{cases}
        \end{equation*}

        \item Temporarily denote
        \begin{equation*}
            D_{X_2}X_2 = uX_1 + vX_2.
        \end{equation*}
        Find $u$ and $v$ from the system of equations:
        \begin{equation*}
        \begin{cases}
            g\left( D_{X_2}X_2,X_1 \right) = B,\\
            g\left(D_{X_2}X_2,X_2 \right) = 0.
        \end{cases}
        \end{equation*}
        \begin{equation*}
        \begin{cases}
            g\left( uX_1 + vX_2,X_1 \right) = B,\\
            g\left(uX_1 + vX_2,X_2 \right) = 0.
        \end{cases}
        \Leftrightarrow
            \begin{cases}
                ua_{11} + va_{12} = B,\\
                ua_{12} + va_{22} = 0,
            \end{cases}
        \Leftrightarrow
        \begin{cases}
            u = \frac{a_{22}B}{\det{(\mathcal{A})}},\\
            v = \frac{-a_{12}B}{\det{(\mathcal{A})}}.
        \end{cases}
        \end{equation*}
    \end{itemize}
    
    Therefore
    \begin{align*}
    & D_{X_1}X_1 = \frac{A}{\det{(\mathcal{A})}}\left( a_{12}X_1 - a_{11}X_2 \right),\\
    & D_{X_1}X_2 = \frac{A}{\det{(\mathcal{A})}}\left( a_{22}X_1 - a_{12}X_2 \right),\\
    & D_{X_2}X_1 = \frac{B}{\det{(\mathcal{A})}}\left( a_{12}X_1 - a_{11}X_2 \right),\\
    & D_{X_2}X_2 = \frac{B}{\det{(\mathcal{A})}}\left( a_{22}X_1 - a_{12}X_2 \right).
    \end{align*}
    \item[$(3)$] \begin{align*}
    & R_{X_1 X_2}X_1 = D_{[X_1,X_2]}X_1 - (D_{X_1}D_{X_2}-D_{X_2}D_{X_1})X_1 = \\  
    & = D_{c_{12}^{1}X_1+c_{12}^{2}X_2}X_1 - D_{X_1}\frac{B}{\det{(\mathcal{A})}}\left( a_{12}X_1 - a_{11}X_2 \right) + D_{X_2}\frac{A}{\det{(\mathcal{A})}}\left( a_{12}X_1 - a_{11}X_2 \right) = \\
    & = D_{c_{12}^{1}X_1+c_{12}^{2}X_2}X_1 - \frac{a_{12}}{\det{(\mathcal{A})}}D_{X_1}\left( BX_1 \right) + \frac{a_{11}}{\det{(\mathcal{A})}}D_{X_1}\left( BX_2 \right) + \frac{a_{12}}{\det{(\mathcal{A})}}D_{X_2}\left( A X_1 \right) - \frac{a_{11}}{\det{(\mathcal{A})}}D_{X_2}\left( A X_2 \right).
    \end{align*}
    
    Let's compute each term separately, and then combine them into one formula.

    \begin{align*}
    & D_{c_{12}^{1}X_1+c_{12}^{2}X_2}X_1 = c_{12}^1D_{X_1}X_1 + c_{12}^2D_{X_2}X_1 = c_{12}^1\frac{A}{\det{(\mathcal{A})}}\left( a_{12}X_1 - a_{11}X_2 \right) + c_{12}^2\frac{B}{\det{(\mathcal{A})}}\left( a_{12}X_1 - a_{11}X_2 \right) = \\
    & = \frac{1}{\det{(\mathcal{A})}}\left( a_{12}\left( c_{12}^1A + c_{12}^2B \right)X_1 - a_{11}\left( c_{12}^1A + c_{12}^2B \right)X_2 \right).
    \end{align*}

    \begin{align*}
    & D_{X_1}\left( B X_1 \right) = D_{X_1}\left( \left( a_{12}c_{12}^1+a_{22}c_{12}^2 \right) X_1 \right) = a_{12}\left[ (X_1c_{12}^1)X_1 + c_{12}^1D_{X_1}X_1 \right] + a_{22}\left[ (X_1 c_{12}^2)X_1 + c_{12}^2D_{X_1}X_1\right] = \\
    & =  a_{12}(X_1c_{12}^1)X_1 + a_{22}(X_1 c_{12}^2)X_1 + a_{12}c_{12}^1\frac{A}{\det{(\mathcal{A})}}\left( a_{12}X_1 - a_{11}X_2 \right) + a_{22}c_{12}^2\frac{A}{\det{(\mathcal{A})}}\left( a_{12}X_1 - a_{11}X_2 \right)=\\
    & = \left[ a_{12}(X_1c_{12}^1) + a_{22}(X_1 c_{12}^2) + \frac{A}{\det{(\mathcal{A})}}(c_{12}^1(a_{12})^2 + c_{12}^2a_{12}a_{22}) \right]X_1 - \frac{A}{\det{(\mathcal{A})}}\left[ c_{12}^1a_{11}a_{12} + c_{12}^2a_{11}a_{22} \right]X_2 = \\
    & = \left[ a_{12}(X_1c_{12}^1) + a_{22}(X_1 c_{12}^2) + \frac{a_{12}AB}{\det{(\mathcal{A})}} \right]X_1 - \frac{a_{11}AB}{\det{(\mathcal{A})}}X_2.
    \end{align*}

    \begin{align*}
    & D_{X_1}\left( B X_2 \right) = D_{X_1}\left( \left( a_{12}c_{12}^1+a_{22}c_{12}^2 \right) X_2 \right) = a_{12}\left[ \left( X_1c_{12}^1 \right)X_2 + c_{12}^1D_{X_1}X_2 \right] + a_{22}\left[ \left( X_1c_{12}^2 \right)X_2 + c_{12}^2D_{X_1}X_2 \right] = \\
    & = a_{12}\left( X_1c_{12}^1 \right)X_2 + a_{22}\left( X_1c_{12}^2 \right)X_2 + a_{12}c_{12}^1\frac{A}{\det{(\mathcal{A})}}\left( a_{22}X_1 - a_{12}X_2 \right) + a_{22}c_{12}^2\frac{A}{\det{(\mathcal{A})}}\left( a_{22}X_1 - a_{12}X_2 \right)=\\
    & = \frac{A}{\det{(\mathcal{A})}}\left[ c_{12}^1a_{12}a_{22} + c_{12}^2\left(a_{22}\right)^2 \right]X_1 + \left[ a_{12}\left( X_1c_{12}^1 \right) + a_{22}\left( X_1c_{12}^2 \right) - \frac{A}{\det{(\mathcal{A})}}\left( c_{12}^1\left(a_{12}\right)^2 + c_{12}^2a_{12}a_{22} \right) \right]X_2 = \\
    & = \frac{a_{22}AB}{\det{(\mathcal{A})}}X_1 + \left[ a_{12}\left( X_1c_{12}^1 \right) + a_{22}\left( X_1c_{12}^2 \right) - \frac{a_{12}AB}{\det{(\mathcal{A})}} \right]X_2.
    \end{align*}

    \begin{align*}
    & D_{X_2}\left( A X_1 \right) = D_{X_2}\left( \left( a_{11}c_{12}^1 + a_{12}c_{12}^2 \right) X_1 \right) = a_{11}\left[ \left( X_2 c_{12}^1 \right)X_1 + c_{12}^1D_{X_2}X_1 \right] + a_{12}\left[ \left( X_2 c_{12}^2 \right)X_1 + c_{12}^2D_{X_2}X_1 \right] = \\
    & = a_{11}\left( X_2 c_{12}^1 \right)X_1 + a_{12}\left( X_2 c_{12}^2 \right)X_1 + a_{11}c_{12}^1\frac{B}{\det{(\mathcal{A})}}\left( a_{12}X_1 - a_{11}X_2 \right) + a_{12}c_{12}^2\frac{B}{\det{(\mathcal{A})}}\left( a_{12}X_1 - a_{11}X_2 \right) = \\
    & = \left[ a_{11}\left( X_2 c_{12}^1 \right) + a_{12}\left( X_2 c_{12}^2 \right) + \frac{B}{\det{(\mathcal{A})}}\left( c_{12}^1a_{11}a_{12} + c_{12}^2\left( a_{12}\right)^2 \right)  \right]X_1 - \frac{B}{\det{(\mathcal{A})}}\left[ c_{12}^1\left( a_{11} \right)^2 + c_{12}^2a_{11}a_{12} \right]X_2 = \\
    & = \left[ a_{11}\left( X_2 c_{12}^1 \right) + a_{12}\left( X_2 c_{12}^2 \right) + \frac{a_{12}AB}{\det{(\mathcal{A})}}  \right]X_1 - \frac{a_{11}AB}{\det{(\mathcal{A})}}X_2.
    \end{align*}
    
    \begin{align*}
    & D_{X_2}\left( A X_2 \right) = D_{X_2}\left( \left( a_{11}c_{12}^1 + a_{12}c_{12}^2 \right)X_2 \right) = a_{11}\left[ \left( X_2c_{12}^1 \right)X_2 + c_{12}^1D_{X_2}X_2 \right] + a_{12}\left[ \left( X_2c_{12}^2 \right)X_2 
 +c_{12}^2D_{X_2}X_2 \right] = \\
 & = a_{11}\left( X_2c_{12}^1 \right)X_2 + a_{12}\left( X_2c_{12}^2 \right)X_2 + a_{11}c_{12}^1\frac{B}{\det{(\mathcal{A})}}\left( a_{22}X_1 - a_{12}X_2 \right) + a_{12}c_{12}^2\frac{B}{\det{(\mathcal{A})}}\left( a_{22}X_1 - a_{12}X_2 \right) = \\
 & = \frac{B}{\det{(\mathcal{A})}}\left[ c_{12}^1a_{11}a_{22} + c_{12}^2a_{12}a_{22} \right]X_1 + \left[ a_{11}\left( X_2c_{12}^1 \right) + a_{12}\left( X_2c_{12}^2 \right) - \frac{B}{\det{(\mathcal{A})}}\left( c_{12}^1a_{11}a_{12} + c_{12}^2 \left(a_{12}\right)^2 \right) \right]X_2 = \\
 & = \frac{a_{22}AB}{\det{(\mathcal{A})}} X_1 + \left[ a_{11}\left( X_2c_{12}^1 \right) + a_{12}\left( X_2c_{12}^2 \right) - \frac{a_{12}AB}{\det{(\mathcal{A})}} \right]X_2.
    \end{align*}

    We obtain

    \begin{align*}
    &  D_{c_{12}^{1}X_1+c_{12}^{2}X_2}X_1 - \frac{a_{12}}{\det{(\mathcal{A})}}D_{X_1}\left( BX_1 \right) + \frac{a_{11}}{\det{(\mathcal{A})}}D_{X_1}\left( BX_2 \right) + \frac{a_{12}}{\det{(\mathcal{A})}}D_{X_2}\left( A X_1 \right) - \frac{a_{11}}{\det{(\mathcal{A})}}D_{X_2}\left( A X_2 \right) = \\
    & = \frac{1}{\det{(\mathcal{A})}}\left( a_{12}\left( c_{12}^1A + c_{12}^2B \right)X_1 - a_{11}\left( c_{12}^1A + c_{12}^2B \right)X_2 \right) - \\
    & - \frac{a_{12}}{\det{(\mathcal{A})}}\left( \left[ a_{12}(X_1c_{12}^1) + a_{22}(X_1 c_{12}^2) + \frac{a_{12}AB}{\det{(\mathcal{A})}} \right]X_1 - \frac{a_{11}AB}{\det{(\mathcal{A})}}X_2 \right)+\\
    & + \frac{a_{11}}{\det{(\mathcal{A})}}\left( \frac{a_{22}AB}{\det{(\mathcal{A})}}X_1 + \left[ a_{12}\left( X_1c_{12}^1 \right) + a_{22}\left( X_1c_{12}^2 \right) - \frac{a_{12}AB}{\det{(\mathcal{A})}} \right]X_2 \right) +\\
    & + \frac{a_{12}}{\det{(\mathcal{A})}}\left( \left[ a_{11}\left( X_2 c_{12}^1 \right) + a_{12}\left( X_2 c_{12}^2 \right) + \frac{a_{12}AB}{\det{(\mathcal{A})}}  \right]X_1 - \frac{a_{11}AB}{\det{(\mathcal{A})}}X_2 \right) - \\
    & - \frac{a_{11}}{\det{(\mathcal{A})}}\left( \frac{a_{22}AB}{\det{(\mathcal{A})}} X_1 + \left[ a_{11}\left( X_2c_{12}^1 \right) + a_{12}\left( X_2c_{12}^2 \right) - \frac{a_{12}AB}{\det{(\mathcal{A})}} \right]X_2 \right) = \\
    & = \frac{1}{\det{\mathcal{A}}}\left(a_{12} \left[ \left( c_{12}^1A + c_{12}^2B \right) - \left( X_1 B \right) + \left( X_2 A \right) \right]X_1 + a_{11} \left[ - \left( c_{12}^1A + c_{12}^2B \right) + \left( X_1B \right) - \left( X_2A \right) \right]X_2 \right) = \frac{1}{\det{\mathcal{A}}}\left( \xi_1 X_1 + \xi_2 X_2 \right).
    \end{align*}
    \item[$(4)$]
    \begin{equation*}
    Q(X_1,X_2) = g(X_1,X_1)g(X_2,X_2)-(g(X_1,X_2))^2 = a_{11} a_{22} - (a_{12})^2 = \det{\mathcal{A}}.
    \end{equation*}
    \item[$(5)$]
    \begin{equation*}
     K = \frac{g(R_{X_1 X_2}X_1,X_2)}{Q(X_1,X_2)} = \frac{g\left( \xi_1 X_1 + \xi_2 X_2,X_2 \right)}{\det{\mathcal{A}}^2} = \frac{a_{12}\xi_1 + a_{22}\xi_2}{\det{\mathcal{A}}^2} = \frac{a_{12}\xi_1 + a_{22}\xi_2}{\det{\mathcal{A}}^2}.
     \end{equation*}
     
\end{itemize}

\section{Three-Dimensional Manifold, Orthonormal Case}

Due to the volume of computations, we managed to derive the sectional curvature formula only for the orthonormal case.

Basis vector fields $X_1$, $X_2$, $X_3$ of the Lorentzian metric $g$.

The metric matrix has the following form:

\begin{equation*}
\mathcal{A} =
\begin{pmatrix}
a_{11} & a_{12} & a_{13} \\
a_{21} & a_{22} & a_{23} \\
a_{31} & a_{32} & a_{33} \\
\end{pmatrix}
=
\begin{pmatrix}
g(X_1,X_1) & g(X_1,X_2) & g(X_1,X_3) \\
g(X_2,X_1) & g(X_2,X_2) & g(X_2,X_3) \\
g(X_3,X_1) & g(X_3,X_2) & g(X_3,X_3) \\
\end{pmatrix}
=
\begin{pmatrix}
-1 & 0 & 0 \\
0 & 1 & 0 \\
0 & 0 & 1
\end{pmatrix}.
\end{equation*}

\subsection{Plan for Computing the Sectional Curvature for a Three-Dimensional Manifold}

\begin{itemize}
    \item[$(1)$] Compute the commutators $[X_i, X_j]$, $i = 1,2,3$, $j = 1,2,3$;
    \item[$(2)$] Compute the Levi-Civita connection coefficients using the Koszul formula (\ref{koszul}), using the basis vector fields $X_1, X_2$, $X_3$;
    \item[$(3)$] Compute the Riemann curvature tensor $R_{v w}v$, where $v$, $w$ define some plane $P$ in the tangent space at point $q$;
    \item[$(4)$] Compute $Q(v,w)$, where $v$, $w$ define some plane $P$ in the tangent space at point $q$;
    \item[$(5)$] Obtain the sectional curvature $K$ using the performed computations and formula (\ref{sectionalcurvature}).
\end{itemize}

\subsection{Main Formula for the Orthonormal Case and Examples}

\begin{theorem}
The sectional curvature of a Lorentzian three-dimensional manifold $M$ with Lorentzian metric $g$, its orthonormal basis given by vector fields $X_1$, $X_2$, $X_3$, $[X_i,X_j] = c_{ij}^{1}X_1+c_{ij}^{2}X_2 + c_{ij}^3X_3$, for an arbitrary plane field defined by the fields $ v = v^1X_1 + v^2X_2 + v^3X_3, \quad w = w^1X_1 + w^2X_2 + w^3X_3$, equals
\begin{equation}
    \label{treh_kriv_ortonorm}
    K = \frac{g\left( R_{v,w}v, w \right)}{Q(v,w)} = \frac{-w^1R_{v,w,v}^1 + w^2R_{v,w,v}^2 + w^3R_{v,w,v}^3}{\left( -\left( v^1 \right)^2 +  \sum_{k=2}^3\left( v^k \right)^2 \right) \cdot \left( -\left( w^1 \right)^2 +  \sum_{l=2}^3\left( w^l \right)^2 \right) - \left( -v^1w^1 + v^2w^2 + v^3w^3 \right)^2},
    \end{equation}

    where

     \begin{equation*}
    R_{v,w,v}^k = \sum_{l, m=1}^3 \left( w^lD_{v,v,l,m}^k - v^lD_{w,v,l,m}^k \right) +\sum_{l=1}^3\left( \sum_{m=1}^2\sum_{m<j\leqslant 3}\left( v^mw^j -v^jw^m \right)c_{m,j}^l + \sum_{m=1}^3 \left( \left( \left( v^mX_m \right)w^l \right) - w^m\left( X_m v^l \right) \right) \right) D_{v,l}^k, \ k = 1,2,3;
    \end{equation*}
    \begin{align*}
    & D_{v,1}^1 = v^2c_{12}^1 + v^3c_{13}^1 + \left( X_1 v^1 \right),\ D_{v,1}^2 = v^1c_{12}^1 - \frac{v^3}{2}\left( c_{32}^1 - c_{21}^3 - c_{13}^2 \right) + \left( X_1 v^2 \right),\ D_{v,1}^3 = v^1c_{13}^1 + \frac{v^2}{2}\left( c_{23}^1 + c_{31}^2 + c_{12}^3 \right) + \left( X_1 v^3 \right);\\
    & D_{v,2}^1 = -v^2c_{12}^2 + \frac{v^3}{2}\left(c_{31}^2 - c_{12}^3 + c_{23}^1 \right) + \left( X_2 v^1 \right),\ D_{v,2}^2 = -v^1c_{12}^2 + v^3c_{23}^2 + \left( X_2 v^2 \right),\ D_{v,2}^3 = -\frac{v^1}{2}\left( c_{13}^2 + c_{32}^1 - c_{21}^3 \right) - v^2c_{23}^2 + \left( X_2 v^3 \right),\\
    & D_{v,3}^1 = \frac{v^2}{2}\left( c_{21}^3-c_{13}^2 + c_{32}^1 \right) - v^3c_{13}^3 + \left( X_3 v^1 \right),\ D_{v,3}^2 = -\frac{v^1}{2}\left( c_{12}^3  + c_{23}^1 - c_{31}^2 \right) + v^3c_{23}^3 + \left( X_3 v^2 \right),\ D_{v,3}^3 = -v^1c_{13}^3 - v^2c_{23}^3 + \left( X_3 v^3 \right); \\
    & D_{w,v,1,1}^1 = w^1D_{v,1}^2c_{12}^1 + w^1D_{v,1}^3c_{13}^1 + \left( X_1 w^1 D_{v,1}^1 \right),\ D_{w,v,1,1}^2 = w^1D_{v,1}^1c_{12}^1 - \frac{w^1D_{v,1}^3}{2}\left( c_{32}^1 - c_{21}^3 - c_{13}^2 \right) + \left( X_1 w^1 D_{v,1}^2  \right),\\ & D_{w,v,1,1}^3 = w^1D_{v,1}^1c_{13}^1 + \frac{w^1D_{v,1}^2}{2}\left( c_{23}^1 +  c_{31}^2 + c_{12}^3  \right) + \left( X_1 w^1 D_{v,1}^3 \right); \\
    & D_{w,v,2,1}^1 = -w^1D_{v,1}^2c_{12}^2 + \frac{w^1D_{v,1}^3}{2}\left(c_{31}^2 - c_{12}^3 + c_{23}^1 \right) + \left( X_2 w^1 D_{v,1}^1 \right),\ D_{w,v,2,1}^2 = -w^1D_{v,1}^1c_{12}^2 + w^1D_{v,1}^3c_{23}^2 + \left( X_2 w^1 D_{v,1}^2  \right),\\ & D_{w,v,2,1}^3 = -\frac{w^1D_{v,1}^1}{2}\left( c_{13}^2 + c_{32}^1 - c_{21}^3 \right) - w^1D_{v,1}^2c_{23}^2 + \left( X_2 w^1 D_{v,1}^3 \right); \\
    & D_{w,v,3,1}^1 = \frac{w^1D_{v,1}^2}{2}\left( c_{21}^3-c_{13}^2 + c_{32}^1 \right) - w^1D_{v,1}^3c_{13}^3 + \left( X_3 w^1 D_{v,1}^1 \right),\ D_{w,v,3,1}^2 = -\frac{w^1D_{v,1}^1}{2}\left( c_{12}^3  + c_{23}^1 - c_{31}^2 \right) + w^1D_{v,1}^3c_{23}^3 + \left( X_3 w^1 D_{v,1}^2  \right),\\ & D_{w,v,3,1}^3 = -w^1D_{v,1}^1c_{13}^3 - w^1D_{v,1}^2c_{23}^3 + \left( X_3 w^1 D_{v,1}^3 \right);\\
    & D_{w,v,1,2}^1 = w^2D_{v,2}^2c_{12}^1 + w^2D_{v,2}^3c_{13}^1 + \left( X_1 w^2 D_{v,2}^1 \right),\ D_{w,v,1,2}^2 = w^2D_{v,2}^1c_{12}^1 - \frac{w^2D_{v,2}^3}{2}\left( c_{32}^1 - c_{21}^3 - c_{13}^2 \right) + \left( X_1 w^2 D_{v,2}^2  \right),\\
    & D_{w,v,1,2}^3 = w^2D_{v,2}^1c_{13}^1 + \frac{w^2D_{v,2}^2}{2}\left( c_{23}^1 + c_{31}^2 + c_{12}^3 \right) + \left( X_1 w^2 D_{v,2}^3 \right);\\
    & D_{w,v,2,2}^1 = -w^2D_{v,2}^2c_{12}^2 + \frac{w^2D_{v,2}^3}{2}\left(c_{31}^2 - c_{12}^3 + c_{23}^1 \right) + \left( X_2 w^2 D_{v,2}^1 \right),\ D_{w,v,2,2}^2 = -w^2D_{v,2}^1c_{12}^2 + w^2D_{v,2}^3c_{23}^2 + \left( X_2 w^2 D_{v,2}^2  \right),\\
    & D_{w,v,2,2}^3 = -\frac{w^2D_{v,2}^1}{2}\left( c_{13}^2 + c_{32}^1 - c_{21}^3 \right) - w^2D_{v,2}^2c_{23}^2 + \left( X_2 w^2 D_{v,2}^3 \right);\\
    & D_{w,v,3,2}^1 = \frac{w^2D_{v,2}^2}{2}\left( c_{21}^3-c_{13}^2 + c_{32}^1 \right) - w^2D_{v,2}^3c_{13}^3 + \left( X_3 w^2 D_{v,2}^1 \right),\ D_{w,v,3,2}^2 = -\frac{w^2D_{v,2}^1}{2}\left( c_{12}^3  + c_{23}^1 - c_{31}^2 \right) + w^2D_{v,2}^3c_{23}^3 + \left( X_3 w^2 D_{v,2}^2  \right),\\
    & D_{w,v,3,2}^3 = -w^2D_{v,2}^1c_{13}^3 - w^2D_{v,2}^2c_{23}^3 + \left( X_3 w^2 D_{v,2}^3 \right);\\
    & D_{w,v,1,3}^1 = w^3D_{v,3}^2c_{12}^1 + w^3D_{v,3}^3c_{13}^1 + \left( X_1 w^3 D_{v,3}^1 \right),\ D_{w,v,1,3}^2 = w^3D_{v,3}^1c_{12}^1 - \frac{w^3D_{v,3}^3}{2}\left( c_{32}^1 - c_{21}^3 - c_{13}^2 \right) + \left( X_1 w^3 D_{v,3}^2  \right),\\
    & D_{w,v,1,3}^3 = w^3D_{v,3}^1c_{13}^1 + \frac{w^3D_{v,3}^2}{2}\left(  c_{23}^1 + c_{31}^2 + c_{12}^3 \right) + \left( X_1 w^3 D_{v,3}^3 \right);\\
    & D_{w,v,2,3}^1 = -w^3D_{v,3}^2c_{12}^2 + \frac{w^3D_{v,3}^3}{2}\left(c_{31}^2 - c_{12}^3 + c_{23}^1 \right) + \left( X_2 w^3 D_{v,3}^1 \right),\ D_{w,v,2,3}^2 = -w^3D_{v,3}^1c_{12}^2 + w^3D_{v,3}^3c_{23}^2 + \left( X_2 w^3 D_{v,3}^2  \right),\\
    & D_{w,v,2,3}^3 = -\frac{w^3D_{v,3}^1}{2}\left( c_{13}^2 + c_{32}^1 - c_{21}^3 \right) - w^3D_{v,3}^2c_{23}^2 + \left( X_2 w^3 D_{v,3}^3 \right);\\
    & D_{w,v,3,3}^1 = \frac{w^3D_{v,3}^2}{2}\left( c_{21}^3-c_{13}^2 + c_{32}^1 \right) - w^3D_{v,3}^3c_{13}^3 + \left( X_3 w^3 D_{v,3}^1 \right),\ D_{w,v,3,3}^2 = -\frac{w^3D_{v,3}^1}{2}\left( c_{12}^3  + c_{23}^1 - c_{31}^2 \right) + w^3D_{v,3}^3c_{23}^3 + \left( X_3 w^3 D_{v,3}^2  \right),\\
    & D_{w,v,3,3}^3 = -w^3D_{v,3}^1c_{13}^3 - w^3D_{v,3}^2c_{23}^3 + \left( X_3 w^3 D_{v,3}^3 \right).
    \end{align*}

    For $D_{v,v,l,m}^k$, $l,m,k = 1,2,3$, the formulas are analogous to those for $D_{w,v,l,m}^k$, only one needs to replace $w$ with $v$.

\end{theorem}

\begin{example}
Consider the plane field defined by $v = X_1$, $w = X_2$, i.e. $v^1 = 1 = w^2$, $v^2 = v^3 = w^1 = w^3 = 0$:

    \begin{align*}
    & K =
    \frac{g\left(R_{X_1,X_2}X_1,X_2 \right)}{-1} = - \left( c_{12}^1c_{12}^1 - c_{12}^2c_{12}^2 - \frac{c_{12}^3}{2}\left( c_{12}^3 + c_{23}^1 - c_{31}^2 \right) - \frac{1}{4}\left( c_{13}^2 + c_{32}^1 - c_{21}^3 \right) \left( c_{32}^1 - c_{21}^3 - c_{13}^2 \right) + c_{13}^1c_{23}^2 + \right. \\
    & + \left. \left( X_2 c_{12}^1 \right) + \left( X_1 c_{12}^2 \right) \right).
    \end{align*}

    Consider the plane field defined by $v = X_1$, $w = X_3$, i.e. $v^1 = 1 = w^3$, $v^2 = v^3 = w^1 = w^2 = 0$:

    \begin{align*}
    & K =
    \frac{g\left(R_{X_1,X_3}X_1,X_3 \right)}{-1} = -\left( c_{13}^1c_{13}^1 - \frac{c_{13}^2}{2}\left( c_{13}^2 + c_{32}^1 - c_{21}^3 \right) - c_{13}^3c_{13}^3 + \frac{1}{4}\left( c_{12}^3  + c_{23}^1 - c_{31}^2 \right)\left( c_{23}^1 + c_{31}^2 + c_{12}^3 \right) - c_{12}^1c_{23}^3 + \right. \\
    & \left. + \left( X_1 c_{13}^3 \right) + \left( X_3 c_{13}^1 \right) \right).
    \end{align*}

    Consider the plane field defined by $v = X_2$, $w = X_3$, i.e. $v^2 = 1 = w^3$, $v^1 = v^3 = w^1 = w^2 = 0$:

    \begin{align*}
    & K =
    \frac{g\left(R_{X_2,X_3}X_2,X_3 \right)}{1} =  \frac{c_{23}^1}{2}\left( c_{23}^1 + c_{31}^2 + c_{12}^3 \right) - c_{23}^2c_{23}^2 - c_{23}^3c_{23}^3 + \frac{1}{4}\left( c_{21}^3-c_{13}^2 + c_{32}^1 \right)\left( c_{13}^2 + c_{32}^1 - c_{21}^3 \right) + c_{12}^2c_{13}^3 + \\
    & + \left( X_2 c_{23}^3 \right) - \left( X_3 c_{23}^2 \right) .
    \end{align*}
\end{example}

\begin{example}
Consider the Lorentzian structure on the Heisenberg group $G\cong \R^3_{x, y, z}$ with the orthonormal frame
$$
X_1 = \frac{\partial}{\partial x } -\frac y2\frac{\partial}{\partial z }, 
\qquad 
X_2 = \frac{\partial}{\partial y } +\frac x2\frac{\partial}{\partial z }, 
\qquad 
X_3 =  \frac{\partial}{\partial z }.
$$

Structural constants $c_{ij}^k = 0$, $i,\ j,\ k = 1,2,3$, except $c_{12}^3 = 1 = -c_{21}^3$.

For fields $v$, $w$ with constant coordinate functions $v^i = const$, $w^j = const$, $i,\ j = 1,2,3$, the sectional curvature equals:

\begin{equation}
\label{kriv_gr_geisenberg}
    K = \frac{-\frac{w^1w^1}{4}\left( 3v^2v^2 - v^3v^3 \right) -\frac{w^2w^2}{4}\left( 3v^1v^1 + v^3v^3 \right) + \frac{w^3w^3}{4}\left( v^1v^1 - v^2v^2 \right) + w^1w^2\frac{3v^1v^2}{2} - w^1w^3\frac{v^1v^3}{2} + w^2w^3\frac{v^2v^3}{2}}{\left( -\left( v^1 \right)^2 +  \sum_{k=2}^3\left( v^k \right)^2 \right) \cdot \left( -\left( w^1 \right)^2 +  \sum_{l=2}^3\left( w^l \right)^2 \right) - \left( -v^1w^1 + v^2w^2 + v^3w^3 \right)^2}.
\end{equation}

\end{example}

\subsection{Derivation of the Main Formula for the Orthonormal Frame}

\begin{itemize}
    \item[$(1)$] Compute the commutators $[X_i, X_j]$, $i = 1,2,3$, $j = 1,2,3$;

Let us write the commutators as follows:
    \begin{equation*}
    [X_i,X_j] = c_{ij}^1 X_1 + c_{ij}^2 X_2 + c_{ij}^3 X_3, \quad i,j = 1,2,3.
    \end{equation*}
    We note right away that $c_{ii}^l = 0$,
    $c_{ij}^l = -c_{ji}^l$, $i,j,l = 1,2,3$, by the properties of the commutator of vector fields.
    
    \item[$(2)$] Compute the Levi-Civita connection coefficients using the Koszul formula (\ref{koszul}), using the basis vector fields $X_1, X_2$, $X_3$;

We note right away that the first three terms of this formula are zero, since $g(X_i,X_j)$ are constant numbers.
 
    \begin{align*}
    & G_{ij}^k = 2g \left(D_{X_i}X_j,X_k \right) = -g\left( X_i,[X_j,X_k] \right) + g\left( X_j,[X_k,X_i] \right) + g\left( X_k,[X_i,X_j] \right) = \\
    & = -g\left(X_i,\sum_{l=1}^{3}c_{jk}^lX_l \right) + g\left( X_j, \sum_{l=1}^{3}c_{ki}^{l}X_l \right) + g\left( X_k, \sum_{l=1}^{3}c_{ij}^{l}X_l \right) = \\
    & = 
    \begin{cases}
    -c_{jk}^{i} + c_{ki}^{j} + c_{ij}^{k}, \quad i,j,k = 2,3;\\
    c_{jk}^{1} + c_{k1}^{j} + c_{1j}^{k}, \quad i = 1;\ j,k = 2,3;\\
    -c_{1k}^{i} - c_{ki}^{1} + c_{i1}^{k}, \quad i = 2,3; \ j = 1;\ k = 2,3;\\
    c_{1k}^{1} - c_{k1}^{1} + c_{11}^{k}, \quad i = j = 1;\ k = 2,3;\\
    -c_{j1}^{i} + c_{1i}^{j} - c_{ij}^{1}, \quad i, j = 2,3;\ k = 1;\\
    c_{j1}^{1} + c_{11}^{j} - c_{1j}^{1}, \quad i = 1,\ j = 2,3; \ k = 1;\\
    -c_{11}^{i} - c_{1i}^{1} - c_{i1}^{1}, \quad i = 2,3;\ j = k = 1;\\
    c_{11}^{1} - c_{11}^{1} - c_{11}^{1}, \quad i = j = k = 1.
    \end{cases}
    =
    \begin{cases}
    -c_{jk}^{i} + c_{ki}^{j} + c_{ij}^{k}, \quad i,j,k = 2,3;\\
    c_{jk}^{1} + c_{k1}^{j} + c_{1j}^{k}, \quad i = 1;\ j,k = 2,3;\\
    -c_{1k}^{i} - c_{ki}^{1} + c_{i1}^{k}, \quad i = 2,3; \ j = 1;\ k = 2,3;\\
    2c_{1k}^{1}, \quad i = j = 1;\ k = 2,3;\\
    -c_{j1}^{i} + c_{1i}^{j} - c_{ij}^{1}, \quad i, j = 2,3;\ k = 1;\\
    -2c_{1j}^{1}, \quad i = 1,\ j = 2,3; \ k = 1;\\
    0, \quad i = 2,3;\ j = k = 1;\\
    0, \quad i = j = k = 1.
    \end{cases}
    \end{align*}

    \begin{equation*}
    D_{X_i}X_j = D_{ij}^1X_1 + D_{ij}^2X_2 + D_{ij}^3X_3.
    \end{equation*}

    We obtain

    \begin{equation*}
    \frac{1}{2}G_{ij}^k = g\left( D_{X_i}X_j,X_k \right) = g\left( \sum_{l=1}D_{ij}^lX_l,X_k \right) = D_{ij}^k a_{kk}.
    \end{equation*}

    Therefore

    \begin{align*}
    & D_{X_1}X_1 = \frac{1}{2}\left( -G_{11}^1X_1 + \sum_{l=2}^3 G_{11}^lX_l \right) = \frac{1}{2}\left( 2c_{12}^1X_2 + 2c_{13}^1X_3 \right) = c_{12}^1X_2 + c_{13}^1X_3;\\
    & D_{X_1}X_2 = \frac{1}{2}\left( -G_{12}^1X_1 + \sum_{l=2}^3 G_{12}^lX_l \right) =  \frac{1}{2}\left( -\left( - 2c_{12}^1 \right)X_1 +  \left( c_{23}^1 + c_{31}^2 + c_{12}^3 \right)X_3 \right) = c_{12}^1X_1 + \frac{1}{2}\left( c_{23}^1 + c_{31}^2 + c_{12}^3 \right)X_3;\\
    & D_{X_1}X_3 = \frac{1}{2}\left( -G_{13}^1X_1 + \sum_{l=2}^3 G_{12}^lX_l \right) = \frac{1}{2}\left( -\left( -2c_{13}^1 \right)X_1 + \left( -c_{32}^1 + c_{21}^3 + c_{13}^2 \right)X_2 \right) = c_{13}^1X_1 - \frac{1}{2}\left( c_{32}^1 - c_{21}^3 - c_{13}^2 \right)X_2; \\
    & D_{X_2}X_1 = \frac{1}{2}\left( -G_{21}^1X_1 + \sum_{l=2}^3 G_{21}^lX_l \right) = \frac{1}{2}\left( \left( -c_{12}^2 - c_{22}^1 + c_{21}^2 \right)X_2 + \left( -c_{13}^2 - c_{32}^1 + c_{21}^3 \right)X_3 \right) = -c_{12}^2X_2 - \frac{1}{2}\left( c_{13}^2 + c_{32}^1 - c_{21}^3 \right)X_3;\\
    & D_{X_2}X_2 = \frac{1}{2}\left( -G_{22}^1X_1 + \sum_{l=2}^3 G_{22}^lX_l \right) = \frac{1}{2}\left( -\left( -c_{21}^2 + c_{12}^2 - c_{22}^1 \right)X_1 + \left( -c_{23}^2 + c_{32}^2 + c_{22}^3 \right)X_3 \right) = -c_{12}^2 X_1 - c_{23}^2 X_3 ;\\
    & D_{X_2}X_3 = \frac{1}{2}\left( -G_{23}^1X_1 + \sum_{l=2}^3 G_{23}^lX_l \right) = \frac{1}{2}\left( -\left( -c_{31}^2 + c_{12}^3 - c_{23}^1 \right)X_1 + \left( -c_{32}^2 + c_{22}^3 + c_{23}^2 \right)X_2 \right) = \frac{1}{2}\left(c_{31}^2 - c_{12}^3 + c_{23}^1 \right) X_1 + c_{23}^2 X_2;\\
    & D_{X_3}X_1 = \frac{1}{2}\left( -G_{31}^1X_1 + \sum_{l=2}^3G_{31}^lX_l \right) = \frac{1}{2}\left( \left( -c_{12}^3  - c_{23}^1 + c_{31}^2 \right)X_2 + \left( -c_{13}^3 - c_{33}^1 + c_{31}^3  \right)X_3 \right) = -\frac{1}{2}\left( c_{12}^3  + c_{23}^1 - c_{31}^2 \right)X_2 - c_{13}^3X_3;\\
    & D_{X_3}X_2 = \frac{1}{2}\left( -G_{32}^1X_1 + \sum_{l=2}^3G_{32}^lX_l \right) = \frac{1}{2}\left( -\left( -c_{21}^3 + c_{13}^2 - c_{32}^1 \right)X_1 + \left( -c_{23}^3 + c_{33}^2 + c_{32}^3  \right)X_3 \right) = \frac{1}{2}\left( c_{21}^3-c_{13}^2 + c_{32}^1 \right)X_1 - c_{23}^3 X_3;\\
    & D_{X_3}X_3 = \frac{1}{2}\left( -G_{33}^1X_1 + \sum_{l=2}^3G_{33}^lX_l \right) = \frac{1}{2}\left( -\left( -c_{31}^3 + c_{13}^3 - c_{33}^1 \right)X_1 + \left( -c_{32}^3 + c_{23}^3 + c_{33}^2 \right)X_2 \right) = -c_{13}^3X_1 + c_{23}^3 X_2.
    \end{align*}
    
    \item[$(3)$] Compute the Riemann curvature tensor $R_{v w}v$, where $v$, $w$ define some plane $P$ in the tangent space at point $q$;

Let us expand the vectors $v$, $w$, which define the plane $P$, in the basis $X_1$, $X_2$, $X_3$:
    \begin{equation*}
    v = v^1X_1 + v^2X_2 + v^3X_3, \quad w = w^1X_1 + w^2X_2 + w^3X_3.
    \end{equation*}

    Compute the curvature tensor $R_{vw}v$.

    Let's do an auxiliary computation:

    \begin{align*}
    & [v^lX_l,w^kX_k] = w^k[v^lX_l,X_k] + \left( \left( v^lX_l \right) w^k \right)X_k = w^k\left( v^l[X_l,X_k]-\left( X_k v^l \right)X_l \right) + \left( \left( v^lX_l \right) w^k \right)X_k = \\
    & = v^lw^k[X_l,X_k] + \left( \left( v^lX_l \right) w^k \right)X_k - w^k\left( X_k v^l \right)X_l;\\
    & \left[ \sum_{l=1}^3v^lX_l, \sum_{k=1}^3w^kX_k \right] = \sum_{l=1}^3\sum_{k=1}^3\left[v^lX_l, w^kX_k\right] = \sum_{l=1}^3\sum_{k=1}^3\left(  v^lw^k[X_l,X_k] + \left( \left( v^lX_l \right) w^k \right)X_k - w^k\left( X_k v^l \right)X_l \right) = \\
    & = \sum_{m=1}^2\sum_{m<k\leqslant 3}\left( v^mw^k -v^kw^m \right)\left[X_m,X_k \right] + \sum_{j=1}^3\left( \sum_{m=1}^3\left( \left( v^mX_m \right)w^j \right) - \sum_{k=1}^3w^k\left( X_k v^j \right) \right)X_j = \\
    & = \sum_{m=1}^2\sum_{m<k\leqslant 3}\left( v^mw^k -v^kw^m \right)\sum_{i=1}^3c_{m,k}^iX_{i} + \sum_{j=1}^3 \left( \sum_{m=1}^3\left( \left( v^mX_m \right)w^j \right) - \sum_{k=1}^3w^k\left( X_k v^j \right) \right) X_j = \sum_{l=1}^3 B^l X_l,
   \end{align*}
    where
    \begin{equation*}
    B^l = \sum_{m=1}^2\sum_{m<j\leqslant 3}\left( v^mw^j -v^jw^m \right)c_{m,j}^l + \sum_{m=1}^3\left( \left( v^mX_m \right)w^l \right) - \sum_{m=1}^3w^m\left( X_m v^l \right), \quad l=1,2,3.
    \end{equation*}

    \begin{align*}
    & R_{vw}v = D_{[v,w]}v-[D_{v},D_{w}]v = D_{\left[ \sum_{l=1}^3v^lX_l, \sum_{k=1}^3w^kX_k \right]}\left(  \sum_{l=1}^3v^lX_l \right) - \left( D_vD_w - D_wD_v \right) \left( \sum_{l=1}^3v^lX_l \right) = \\
    & = D_{\sum_{m=1}^3B^m X_m}\left( \sum_{l=1}^3v^lX_l  \right) - D_vD_w\left( \sum_{l=1}^3v^lX_l \right) + D_wD_v\left( \sum_{l=1}^3v^lX_l \right) = \\
    & = \left( \sum_{j=1}^3B^jD_{X_j} \right) \left( \sum_{l=1}^3v^lX_l \right) - D_v\left( \sum_{l=1}^3w^lD_{X_l} \right)\left( \sum_{l=1}^3v^lX_l \right) + D_w\left( \sum_{l=1}^3v^lD_{X_l} \right) \left( \sum_{l=1}^3v^lX_l \right) = \\
    & = \left( \sum_{j=1}^3B^jD_{X_j} \right) \left( \sum_{l=1}^3v^lX_l \right) - \left( \sum_{l=1}^3v^lD_{X_l} \right)\left( \sum_{m=1}^3w^mD_{X_m} \right)\left( \sum_{k=1}^3v^kX_k \right) + \left( \sum_{l=1}^3w^lD_{X_l} \right)\left( \sum_{m=1}^3v^mD_{X_m} \right) \left( \sum_{k=1}^3v^kX_k \right) = \\
    & = \sum_{l=1}^{3}B^l \sum_{k=1}^3 D_{v,l}^kX_k - \sum_{l=1}^3v^l\sum_{m, k=1}^3 D_{w,v,l,m}^kX_k + \sum_{l=1}^3w^l\sum_{m, k=1}^3 D_{v,v,l,m}^kX_k = \\
    & = \left( \sum_{l,m=1}^3\left( w^lD_{v,v,l,m}^1 - v^lD_{w,v,l,m}^1 \right) + \sum_{l=1}^{3}B^l D_{v,l}^1 \right) X_1 + \left( \sum_{l,m=1}^3\left( w^lD_{v,v,l,m}^2 - v^lD_{w,v,l,m}^2 \right) + \sum_{l=1}^{3}B^l D_{v,l}^2 \right) X_2 + \\
    & + \left( \sum_{l,m=1}^3\left( w^lD_{v,v,l,m}^3 - v^lD_{w,v,l,m}^3 \right) + \sum_{l=1}^{3}B^l D_{v,l}^3 \right) X_3 = R_{v,w,v}^1X_1 + R_{v,w,v}^2X_2 + R_{v,w,v}^3X_3.  
    \end{align*}

    For the general case, additional computations will be required:

    \begin{align*}
    & D_{X_i}\left( v \right) = D_{X_i}\left( v^1X_1 + v^2X_2 + v^3X_3 \right) = \left( X_i v^1 \right) X_1 + \left( X_i v^2 \right) X_2 + \left( X_i v^3 \right) X_3 + v^1 D_{X_i}X_1 + v^2 D_{X_i}X_2 + v^3 D_{X_i}X_3 = \\
    & = D_{v,i}^1X_1 + D_{v,i}^2X_2 + D_{v,i}^3X_3.
    \end{align*}

    \begin{align*}
    & D_{X_i}\left( w^jD_{X_j}v \right) = D_{X_i}\left( w^j\left( D_{v,j}^1X_1 + D_{v,j}^2X_2 + D_{v,j}^3X_3 \right) \right) = \left( X_i w^j D_{v,j}^1 \right)X_1 +  \left( X_i w^j D_{v,j}^2  \right)X_2 +  \left( X_i w^j D_{v,j}^3 \right)X_3 +  \\
    & + w^jD_{v,j}^1 D_{X_i}X_1 +  w^jD_{v,j}^2D_{X_i}X_2 + w^jD_{v,j}^3 D_{X_i}X_3 = D_{w,v,i,j}^1X_1 + D_{w,v,i,j}^2X_2 + D_{w,v,i,j}^3X_3. 
    \end{align*}

    \begin{align*}
    & D_{X_i}\left( v^jD_{X_j}v \right) = D_{X_i}\left( v^j\left( D_{v,j}^1X_1 + D_{v,j}^2X_2 + D_{v,j}^3X_3 \right) \right) = \left( X_i v^j D_{v,j}^1 \right)X_1 +  \left( X_i v^j D_{v,j}^2  \right)X_2 +  \left( X_i v^j D_{v,j}^3 \right)X_3 + \\
    & + v^jD_{v,j}^1 D_{X_i}X_1 +  v^jD_{v,j}^2D_{X_i}X_2 + v^jD_{v,j}^3 D_{X_i}X_3 = D_{v,v,i,j}^1X_1 + D_{v,v,i,j}^2X_2 + D_{v,v,i,j}^3X_3. 
    \end{align*}

    The detailed expression for $D_{v,l}^k$, $D_{w,v,l,m}^k$, $D_{v,v,l,m}^k$, $l,\ m,\ k = 1,2,3$, in terms of structural and coordinate functions, as well as their derivatives, can be found in Appendix 1.
    
    \item[$(4)$] Compute $Q(v,w)$, where $v$, $w$ define some plane $P$ in the tangent space at point $q$;

     \begin{align*}
    & Q(v, w) = g(v,v)g(w,w) - g(v,w)^2 = \left( -\left( v^1 \right)^2 + \left( v^2 \right)^2 + \left( v^3 \right)^2 \right)\left( -\left( w^1 \right)^2 + \left( w^2 \right)^2 + \left( w^3 \right)^2 \right) - \left( -v^1w^1 + v^2w^2 + v^3w^3 \right)^2 = \\
    & = \left( -\left( v^1 \right)^2 +  \sum_{k=2}^3\left( v^k \right)^2 \right) \cdot \left( -\left( w^1 \right)^2 +  \sum_{l=2}^3\left( w^l \right)^2 \right) - \left( -v^1w^1 + v^2w^2 + v^3w^3 \right)^2.
    \end{align*}
    
    \item[$(5)$] Obtain the sectional curvature $K$ using the performed computations and formula (\ref{sectionalcurvature}).

     \begin{align*}
    & K = \frac{g\left( R_{v,w}v, w \right)}{Q(v,w)} = \frac{-w^1R_{v,w,v}^1 + w^2R_{v,w,v}^2 + w^3R_{v,w,v}^3}{\left( -\left( v^1 \right)^2 +  \sum_{k=2}^3\left( v^k \right)^2 \right) \cdot \left( -\left( w^1 \right)^2 +  \sum_{l=2}^3\left( w^l \right)^2 \right) - \left( -v^1w^1 + v^2w^2 + v^3w^3 \right)^2},
    \end{align*}

    where

    \begin{equation*}
    R_{v,w,v}^k = \sum_{l,m=1}^3\left( w^lD_{v,v,l,m}^k - v^lD_{w,v,l,m}^k \right) + \sum_{l=1}^{3} B^l  D_{v,l}^k, \quad k = 1,2,3,
    \end{equation*}

    \begin{equation*}
    B^l = \sum_{m=1}^2\sum_{m<j\leqslant 3}\left( v^mw^j -v^jw^m \right)c_{m,j}^l + \sum_{m=1}^3\left( \left( v^mX_m \right)w^l \right) - \sum_{m=1}^3w^m\left( X_m v^l \right), \quad l=1,2,3.    
    \end{equation*}
    
\end{itemize}

\section{Appendix 1: Computation of the Coefficients $D_{v,l}^k$, $D_{w,v,l,m}^k$, $D_{v,v,l,m}^k$}

    \begin{align*}
    & D_{X_1}\left( v \right) =  D_{X_1}\left( v^1X_1 + v^2X_2 + v^3X_3 \right) = \left( X_1 v^1 \right) X_1 + \left( X_1 v^2 \right) X_2 + \left( X_1 v^3 \right) X_3 + v^1 D_{X_1}X_1 + v^2 D_{X_1}X_2 + v^3 D_{X_1}X_3 = \\
    & = \left( X_1 v^1 \right) X_1 + \left( X_1 v^2 \right) X_2 + \left( X_1 v^3 \right) X_3 + v^1 \left( c_{12}^1X_2 + c_{13}^1X_3 \right) + v^2 \left( c_{12}^1X_1 + \frac{1}{2}\left( c_{23}^1 + c_{31}^2 + c_{12}^3 \right)X_3 \right) + \\
    & + v^3 \left( c_{13}^1X_1 - \frac{1}{2}\left( c_{32}^1 - c_{21}^3 - c_{13}^2 \right)X_2 \right) = \left( v^2c_{12}^1 + v^3c_{13}^1 + \left( X_1 v^1 \right) \right)X_1 + \left( v^1c_{12}^1 - \frac{v^3}{2}\left( c_{32}^1 - c_{21}^3 - c_{13}^2 \right) + \left( X_1 v^2 \right) \right)X_2 + \\
    & + \left( v^1c_{13}^1 + \frac{v^2}{2}\left( c_{23}^1 + c_{31}^2 + c_{12}^3 \right) + \left( X_1 v^3 \right) \right)X_3 = D_{v,1}^1X_1 + D_{v,1}^2X_2 + D_{v,1}^3X_3,
    \end{align*}

    \begin{align*}
    & D_{X_2}\left( v \right) = D_{X_2}\left( v^1X_1 + v^2X_2 + v^3X_3 \right) = \left( X_2 v^1 \right) X_1 + \left( X_2 v^2 \right) X_2 + \left( X_2 v^3 \right) X_3 + v^1 D_{X_2}X_1 + v^2 D_{X_2}X_2 + v^3 D_{X_2}X_3 = \\
    & = \left( X_2 v^1 \right) X_1 + \left( X_2 v^2 \right) X_2 + \left( X_2 v^3 \right) X_3 + v^1\left( -c_{12}^2X_2 - \frac{1}{2}\left( c_{13}^2 + c_{32}^1 - c_{21}^3 \right)X_3 \right) + v^2\left( -c_{12}^2 X_1 - c_{23}^2 X_3 \right) +\\
    & + v^3\left( \frac{1}{2}\left(c_{31}^2 - c_{12}^3 + c_{23}^1 \right) X_1 + c_{23}^2 X_2  \right) = \left( -v^2c_{12}^2 + \frac{v^3}{2}\left(c_{31}^2 - c_{12}^3 + c_{23}^1 \right) + \left( X_2 v^1 \right) \right)X_1 + \left( -v^1c_{12}^2 + v^3c_{23}^2 + \left( X_2 v^2 \right) \right)X_2 + \\
    & + \left( -\frac{v^1}{2}\left( c_{13}^2 + c_{32}^1 - c_{21}^3 \right) - v^2c_{23}^2 + \left( X_2 v^3 \right) \right)X_3 = D_{v,2}^1X_1 + D_{v,2}^2X_2 + D_{v,2}^3X_3,
    \end{align*}

    \begin{align*}
    & D_{X_3}\left( v \right) = D_{X_3}\left( v^1X_1 + v^2X_2 + v^3X_3 \right) = \left( X_3 v^1 \right) X_1 + \left( X_3 v^2 \right) X_2 + \left( X_3 v^3 \right) X_3 + v^1 D_{X_3}X_1 + v^2 D_{X_3}X_2 + v^3 D_{X_3}X_3 = \\
    & = \left( X_3 v^1 \right) X_1 + \left( X_3 v^2 \right) X_2 + \left( X_3 v^3 \right) X_3 + v^1 \left( -\frac{1}{2}\left( c_{12}^3  + c_{23}^1 - c_{31}^2 \right)X_2 - c_{13}^3X_3 \right) + v^2 \left( \frac{1}{2}\left( c_{21}^3-c_{13}^2 + c_{32}^1 \right)X_1 - c_{23}^3 X_3 \right) + \\
    & + v^3 \left( -c_{13}^3X_1 + c_{23}^3 X_2 \right) = \left( \frac{v^2}{2}\left( c_{21}^3-c_{13}^2 + c_{32}^1 \right) - v^3c_{13}^3 + \left( X_3 v^1 \right) \right)X_1 + \left( -\frac{v^1}{2}\left( c_{12}^3  + c_{23}^1 - c_{31}^2 \right) + v^3c_{23}^3 + \left( X_3 v^2 \right) \right)X_2 + \\
    & + \left( -v^1c_{13}^3 - v^2c_{23}^3 + \left( X_3 v^3 \right) \right)X_3 = D_{v,3}^1X_1 + D_{v,3}^2X_2 + D_{v,3}^3X_3.
    \end{align*}

    \begin{align*}
    & D_{X_1}\left( w^1D_{X_1}v \right) = D_{X_1}\left( w^1\left( D_{v,1}^1X_1 + D_{v,1}^2X_2 + D_{v,1}^3X_3 \right) \right) = \left( X_1 w^1 D_{v,1}^1 \right)X_1 +  \left( X_1 w^1 D_{v,1}^2  \right)X_2 +  \left( X_1 w^1 D_{v,1}^3 \right)X_3 + \\ & + w^1D_{v,1}^1 D_{X_1}X_1 +  w^1D_{v,1}^2D_{X_1}X_2 + w^1D_{v,1}^3 D_{X_1}X_3 = \left( X_1 w^1 D_{v,1}^1 \right)X_1 +  \left( X_1 w^1 D_{v,1}^2  \right)X_2 +  \left( X_1 w^1 D_{v,1}^3 \right)X_3 + \\
    & + w^1D_{v,1}^1 \left( c_{12}^1X_2 + c_{13}^1X_3 \right) + w^1D_{v,1}^2\left( c_{12}^1X_1 + \frac{1}{2}\left( c_{23}^1 + c_{31}^2 + c_{12}^3 \right)X_3 \right) + w^1D_{v,1}^3\left( c_{13}^1X_1 - \frac{1}{2}\left( c_{32}^1 - c_{21}^3 - c_{13}^2 \right)X_2 \right) = \\
    & = \left( w^1D_{v,1}^2c_{12}^1 + w^1D_{v,1}^3c_{13}^1 + \left( X_1 w^1 D_{v,1}^1 \right) \right) X_1 + \left( w^1D_{v,1}^1c_{12}^1 - \frac{w^1D_{v,1}^3}{2}\left( c_{32}^1 - c_{21}^3 - c_{13}^2 \right) + \left( X_1 w^1 D_{v,1}^2  \right) \right)X_2 + \\
    & + \left( w^1D_{v,1}^1c_{13}^1 + \frac{w^1D_{v,1}^2}{2}\left( c_{23}^1 + c_{31}^2 + c_{12}^3 \right) + \left( X_1 w^1 D_{v,1}^3 \right) \right)X_3 = D_{w,v,1,1}^1X_1 + D_{w,v,1,1}^2X_2 + D_{w,v,1,1}^3X_3,
    \end{align*}

    \begin{align*}
    & D_{X_2}\left( w^1D_{X_1}v \right) = D_{X_2}\left( w^1\left( D_{v,1}^1X_1 + D_{v,1}^2X_2 + D_{v,1}^3X_3 \right) \right) = \left( X_2 w^1 D_{v,1}^1 \right)X_1 +  \left( X_2 w^1 D_{v,1}^2  \right)X_2 +  \left( X_2 w^1 D_{v,1}^3 \right)X_3 + \\ & +  w^1D_{v,1}^1 D_{X_2}X_1 +  w^1D_{v,1}^2D_{X_2}X_2 + w^1D_{v,1}^3 D_{X_2}X_3 = \left( X_2 w^1 D_{v,1}^1 \right)X_1 +  \left( X_2 w^1 D_{v,1}^2  \right)X_2 +  \left( X_2 w^1 D_{v,1}^3 \right)X_3 + \\
    & + w^1D_{v,1}^1 \left( -c_{12}^2X_2 - \frac{1}{2}\left( c_{13}^2 + c_{32}^1 - c_{21}^3 \right)X_3 \right) +  w^1D_{v,1}^2 \left( -c_{12}^2 X_1 - c_{23}^2 X_3 \right) + w^1D_{v,1}^3 \left( \frac{1}{2}\left(c_{31}^2 - c_{12}^3 + c_{23}^1 \right) X_1 + c_{23}^2 X_2 \right) = \\
    & = \left( -w^1D_{v,1}^2c_{12}^2 + \frac{w^1D_{v,1}^3}{2}\left(c_{31}^2 - c_{12}^3 + c_{23}^1 \right) + \left( X_2 w^1 D_{v,1}^1 \right) \right)X_1 + \left( -w^1D_{v,1}^1c_{12}^2 + w^1D_{v,1}^3c_{23}^2 + \left( X_2 w^1 D_{v,1}^2  \right) \right)X_2 + \\
    & + \left( -\frac{w^1D_{v,1}^1}{2}\left( c_{13}^2 + c_{32}^1 - c_{21}^3 \right) - w^1D_{v,1}^2c_{23}^2 + \left( X_2 w^1 D_{v,1}^3 \right) \right)X_3 = D_{w,v,2,1}^1X_1 + D_{w,v,2,1}^2X_2 + D_{w,v,2,1}^3X_3,
    \end{align*}

    \begin{align*}
    & D_{X_3}\left( w^1D_{X_1}v \right) = D_{X_3}\left( w^1\left( D_{v,1}^1X_1 + D_{v,1}^2X_2 + D_{v,1}^3X_3 \right) \right) = \left( X_3 w^1 D_{v,1}^1 \right)X_1 +  \left( X_3 w^1 D_{v,1}^2  \right)X_2 +  \left( X_3 w^1 D_{v,1}^3 \right)X_3 + \\ & +  w^1D_{v,1}^1 D_{X_3}X_1 +  w^1D_{v,1}^2D_{X_3}X_2 + w^1D_{v,1}^3 D_{X_3}X_3 = \left( X_3 w^1 D_{v,1}^1 \right)X_1 +  \left( X_3 w^1 D_{v,1}^2  \right)X_2 +  \left( X_3 w^1 D_{v,1}^3 \right)X_3 + \\ 
    & + w^1D_{v,1}^1 \left( -\frac{1}{2}\left( c_{12}^3  + c_{23}^1 - c_{31}^2 \right)X_2 - c_{13}^3X_3 \right) + 
    w^1D_{v,1}^2 \left( \frac{1}{2}\left( c_{21}^3-c_{13}^2 + c_{32}^1 \right)X_1 - c_{23}^3 X_3 \right) +
    w^1D_{v,1}^3 \left( -c_{13}^3X_1 + c_{23}^3 X_2 \right) = \\
    & = \left( \frac{w^1D_{v,1}^2}{2}\left( c_{21}^3-c_{13}^2 + c_{32}^1 \right) - w^1D_{v,1}^3c_{13}^3 + \left( X_3 w^1 D_{v,1}^1 \right) \right)X_1 +
    \left( -\frac{w^1D_{v,1}^1}{2}\left( c_{12}^3  + c_{23}^1 - c_{31}^2 \right) + w^1D_{v,1}^3c_{23}^3 + \left( X_3 w^1 D_{v,1}^2  \right) \right)X_2 + \\
    & + \left( -w^1D_{v,1}^1c_{13}^3 - w^1D_{v,1}^2c_{23}^3 + \left( X_3 w^1 D_{v,1}^3 \right) \right)X_3 = D_{w,v,3,1}^1X_1 + D_{w,v,3,1}^2X_2 + D_{w,v,3,1}^3X_3.
    \end{align*}

    \begin{align*}
    & D_{X_1}\left( w^2D_{X_2}v \right) = D_{X_1}\left( w^2 \left( D_{v,2}^1X_1 + D_{v,2}^2X_2 + D_{v,2}^3X_3 \right) \right) = \left( X_1 w^2 D_{v,2}^1 \right)X_1 +  \left( X_1 w^2 D_{v,2}^2  \right)X_2 +  \left( X_1 w^2 D_{v,2}^3 \right)X_3 + \\ & + w^2D_{v,2}^1 D_{X_1}X_1 +  w^2D_{v,2}^2D_{X_1}X_2 + w^2D_{v,2}^3 D_{X_1}X_3 = \left( X_1 w^2 D_{v,2}^1 \right)X_1 +  \left( X_1 w^2 D_{v,2}^2  \right)X_2 +  \left( X_1 w^2 D_{v,2}^3 \right)X_3 + \\ 
    & + w^2D_{v,2}^1 \left( c_{12}^1X_2 + c_{13}^1X_3 \right) +  w^2D_{v,2}^2 \left( c_{12}^1X_1 + \frac{1}{2}\left( c_{23}^1 + c_{31}^2 + c_{12}^3 \right)X_3 \right) + w^2D_{v,2}^3 \left( c_{13}^1X_1 - \frac{1}{2}\left( c_{32}^1 - c_{21}^3 - c_{13}^2 \right)X_2 \right) = \\
    & = \left( w^2D_{v,2}^2c_{12}^1 + w^2D_{v,2}^3c_{13}^1 + \left( X_1 w^2 D_{v,2}^1 \right) \right) X_1 + \left( w^2D_{v,2}^1c_{12}^1 - \frac{w^2D_{v,2}^3}{2}\left( c_{32}^1 - c_{21}^3 - c_{13}^2 \right) + \left( X_1 w^2 D_{v,2}^2  \right) \right)X_2 + \\
    & + \left( w^2D_{v,2}^1c_{13}^1 + \frac{w^2D_{v,2}^2}{2}\left( c_{23}^1 + c_{31}^2 + c_{12}^3 \right) + \left( X_1 w^2 D_{v,2}^3 \right) \right)X_3 = D_{w,v,1,2}^1X_1 + D_{w,v,1,2}^2X_2 + D_{w,v,1,2}^3X_3,
    \end{align*}

    \begin{align*}
    & D_{X_2}\left( w^2D_{X_2}v \right) = D_{X_2}\left( w^2\left( D_{v,2}^1X_1 + D_{v,2}^2X_2 + D_{v,2}^3X_3 \right) \right) = \left( X_2 w^2 D_{v,2}^1 \right)X_1 +  \left( X_2 w^2 D_{v,2}^2  \right)X_2 +  \left( X_2 w^2 D_{v,2}^3 \right)X_3 + \\ & +  w^2D_{v,2}^1 D_{X_2}X_1 +  w^2D_{v,2}^2D_{X_2}X_2 + w^2D_{v,2}^3 D_{X_2}X_3 = \left( X_2 w^2 D_{v,2}^1 \right)X_1 +  \left( X_2 w^2 D_{v,2}^2  \right)X_2 +  \left( X_2 w^2 D_{v,2}^3 \right)X_3 + \\ & +  w^2D_{v,2}^1 \left( -c_{12}^2X_2 - \frac{1}{2}\left( c_{13}^2 + c_{32}^1 - c_{21}^3 \right)X_3 \right) +  w^2D_{v,2}^2 \left( -c_{12}^2 X_1 - c_{23}^2 X_3 \right) + w^2D_{v,2}^3 \left( \frac{1}{2}\left(c_{31}^2 - c_{12}^3 + c_{23}^1 \right) X_1 + c_{23}^2 X_2 \right) = \\
    & = \left( -w^2D_{v,2}^2c_{12}^2 + \frac{w^2D_{v,2}^3}{2}\left(c_{31}^2 - c_{12}^3 + c_{23}^1 \right) + \left( X_2 w^2 D_{v,2}^1 \right) \right)X_1 + \left( -w^2D_{v,2}^1c_{12}^2 + w^2D_{v,2}^3c_{23}^2 + \left( X_2 w^2 D_{v,2}^2  \right) \right)X_2 + \\
    & + \left( -\frac{w^2D_{v,2}^1}{2}\left( c_{13}^2 + c_{32}^1 - c_{21}^3 \right) - w^2D_{v,2}^2c_{23}^2 + \left( X_2 w^2 D_{v,2}^3 \right) \right)X_3 = D_{w,v,2,2}^1X_1 + D_{w,v,2,2}^2X_2 + D_{w,v,2,2}^3X_3,
    \end{align*}

    \begin{align*}
    & D_{X_3}\left( w^2D_{X_2}v \right) = D_{X_3}\left( w^2\left( D_{v,2}^1X_1 + D_{v,2}^2X_2 + D_{v,2}^3X_3 \right) \right) = \left( X_3 w^2 D_{v,2}^1 \right)X_1 +  \left( X_3 w^2 D_{v,2}^2  \right)X_2 +  \left( X_3 w^2 D_{v,2}^3 \right)X_3 + \\ & +  w^2D_{v,2}^1 D_{X_3}X_1 +  w^2D_{v,2}^2D_{X_3}X_2 + w^2D_{v,2}^3 D_{X_3}X_3 = \left( X_3 w^2 D_{v,2}^1 \right)X_1 +  \left( X_3 w^2 D_{v,2}^2  \right)X_2 +  \left( X_3 w^2 D_{v,2}^3 \right)X_3 + \\ & +  w^2D_{v,2}^1 \left( -\frac{1}{2}\left( c_{12}^3  + c_{23}^1 - c_{31}^2 \right)X_2 - c_{13}^3X_3 \right) +  w^2D_{v,2}^2 \left( \frac{1}{2}\left( c_{21}^3-c_{13}^2 + c_{32}^1 \right)X_1 - c_{23}^3 X_3 \right) + w^2D_{v,2}^3 \left( -c_{13}^3X_1 + c_{23}^3 X_2 \right) = \\
    & = \left( \frac{w^2D_{v,2}^2}{2}\left( c_{21}^3-c_{13}^2 + c_{32}^1 \right) - w^2D_{v,2}^3c_{13}^3 + \left( X_3 w^2 D_{v,2}^1 \right) \right)X_1 +
    \left( -\frac{w^2D_{v,2}^1}{2}\left( c_{12}^3  + c_{23}^1 - c_{31}^2 \right) + w^2D_{v,2}^3c_{23}^3 + \left( X_3 w^2 D_{v,2}^2  \right) \right)X_2 + \\
    & + \left( -w^2D_{v,2}^1c_{13}^3 - w^2D_{v,2}^2c_{23}^3 + \left( X_3 w^2 D_{v,2}^3 \right) \right)X_3 = D_{w,v,3,2}^1X_1 + D_{w,v,3,2}^2X_2 + D_{w,v,3,2}^3X_3.
    \end{align*}

    \begin{align*}
    & D_{X_1}\left( w^3D_{X_3}v \right) = D_{X_1}\left( w^3 \left( D_{v,3}^1X_1 + D_{v,3}^2X_2 + D_{v,3}^3X_3 \right) \right) = \left( X_1 w^3 D_{v,3}^1 \right)X_1 +  \left( X_1 w^3 D_{v,3}^2  \right)X_2 +  \left( X_1 w^3 D_{v,3}^3 \right)X_3 + \\ 
    & + w^3D_{v,3}^1D_{X_1}X_1 + w^3D_{v,3}^2D_{X_1}X_2 + w^3D_{v,3}^3D_{X_1}X_3 = \left( X_1 w^3 D_{v,3}^1 \right)X_1 +  \left( X_1 w^3 D_{v,3}^2  \right)X_2 +  \left( X_1 w^3 D_{v,3}^3 \right)X_3 + \\ 
    & + w^3D_{v,3}^1 \left( c_{12}^1X_2 + c_{13}^1X_3 \right) +  w^3D_{v,3}^2 \left( c_{12}^1X_1 + \frac{1}{2}\left( c_{23}^1 + c_{31}^2 + c_{12}^3 \right)X_3 \right) + w^3D_{v,3}^3 \left( c_{13}^1X_1 - \frac{1}{2}\left( c_{32}^1 - c_{21}^3 - c_{13}^2 \right)X_2 \right) = \\
    & = \left( w^3D_{v,3}^2c_{12}^1 + w^3D_{v,3}^3c_{13}^1 + \left( X_1 w^3 D_{v,3}^1 \right) \right) X_1 + \left( w^3D_{v,3}^1c_{12}^1 - \frac{w^3D_{v,3}^3}{2}\left( c_{32}^1 - c_{21}^3 - c_{13}^2 \right) + \left( X_1 w^3 D_{v,3}^2  \right) \right)X_2 + \\
    & + \left( w^3D_{v,3}^1c_{13}^1 + \frac{w^3D_{v,3}^2}{2}\left(  c_{23}^1 + c_{31}^2 + c_{12}^3 \right) + \left( X_1 w^3 D_{v,3}^3 \right) \right)X_3 = D_{w,v,1,3}^1X_1 + D_{w,v,1,3}^2X_2 + D_{w,v,1,3}^3X_3,
    \end{align*}

    \begin{align*}
    & D_{X_2}\left( w^3D_{X_3}v \right) = D_{X_2}\left( w^3\left( D_{v,3}^1X_1 + D_{v,3}^2X_2 + D_{v,3}^3X_3 \right) \right) = \left( X_2 w^3 D_{v,3}^1 \right)X_1 +  \left( X_2 w^3 D_{v,3}^2  \right)X_2 +  \left( X_2 w^3 D_{v,3}^3 \right)X_3 + \\ & +  w^3D_{v,3}^1 D_{X_2}X_1 +  w^3D_{v,3}^2D_{X_2}X_2 + w^3D_{v,3}^3 D_{X_2}X_3 = \left( X_2 w^3 D_{v,3}^1 \right)X_1 +  \left( X_2 w^3 D_{v,3}^2  \right)X_2 +  \left( X_2 w^3 D_{v,3}^3 \right)X_3 + \\ & +  w^3D_{v,3}^1 \left( -c_{12}^2X_2 - \frac{1}{2}\left( c_{13}^2 + c_{32}^1 - c_{21}^3 \right)X_3 \right) +  w^3D_{v,3}^2 \left( -c_{12}^2 X_1 - c_{23}^2 X_3 \right) + w^3D_{v,3}^3 \left( \frac{1}{2}\left(c_{31}^2 - c_{12}^3 + c_{23}^1 \right) X_1 + c_{23}^2 X_2 \right) = \\
    & = \left( -w^3D_{v,3}^2c_{12}^2 + \frac{w^3D_{v,3}^3}{2}\left(c_{31}^2 - c_{12}^3 + c_{23}^1 \right) + \left( X_2 w^3 D_{v,3}^1 \right) \right)X_1 + \left( -w^3D_{v,3}^1c_{12}^2 + w^3D_{v,3}^3c_{23}^2 + \left( X_2 w^3 D_{v,3}^2  \right) \right)X_2 + \\
    & + \left( -\frac{w^3D_{v,3}^1}{2}\left( c_{13}^2 + c_{32}^1 - c_{21}^3 \right) - w^3D_{v,3}^2c_{23}^2 + \left( X_2 w^3 D_{v,3}^3 \right) \right)X_3 = D_{w,v,2,3}^1X_1 + D_{w,v,2,3}^2X_2 + D_{w,v,2,3}^3X_3,
    \end{align*}

    \begin{align*}
    & D_{X_3}\left( w^3D_{X_3}v \right) = D_{X_3}\left( w^3\left( D_{v,3}^1X_1 + D_{v,3}^2X_2 + D_{v,3}^3X_3 \right) \right) = \left( X_3 w^3 D_{v,3}^1 \right)X_1 +  \left( X_3 w^3 D_{v,3}^2  \right)X_2 +  \left( X_3 w^3 D_{v,3}^3 \right)X_3 + \\ & +  w^3D_{v,3}^1 D_{X_3}X_1 +  w^3D_{v,3}^2D_{X_3}X_2 + w^3D_{v,3}^3 D_{X_3}X_3 = \left( X_3 w^3 D_{v,3}^1 \right)X_1 +  \left( X_3 w^3 D_{v,3}^2  \right)X_2 +  \left( X_3 w^3 D_{v,3}^3 \right)X_3 + \\ & +  w^3D_{v,3}^1 \left( -\frac{1}{2}\left( c_{12}^3  + c_{23}^1 - c_{31}^2 \right)X_2 - c_{13}^3X_3 \right) +  w^3D_{v,3}^2\left( \frac{1}{2}\left( c_{21}^3-c_{13}^2 + c_{32}^1 \right)X_1 - c_{23}^3 X_3 \right) + w^3D_{v,3}^3 \left( -c_{13}^3X_1 + c_{23}^3 X_2 \right) = \\
    & = \left( \frac{w^3D_{v,3}^2}{2}\left( c_{21}^3-c_{13}^2 + c_{32}^1 \right) - w^3D_{v,3}^3c_{13}^3 + \left( X_3 w^3 D_{v,3}^1 \right) \right)X_1 +
    \left( -\frac{w^3D_{v,3}^1}{2}\left( c_{12}^3  + c_{23}^1 - c_{31}^2 \right) + w^3D_{v,3}^3c_{23}^3 + \left( X_3 w^3 D_{v,3}^2  \right) \right)X_2 + \\
    & + \left( -w^3D_{v,3}^1c_{13}^3 - w^3D_{v,3}^2c_{23}^3 + \left( X_3 w^3 D_{v,3}^3 \right) \right)X_3 = D_{w,v,3,3}^1X_1 + D_{w,v,3,3}^2X_2 + D_{w,v,3,3}^3X_3.
    \end{align*}

    For $D_{v,v,l,m}^k$, the calculations are similar, only instead of the coordinates $w$, we have the coordinates $v$.

\section{Appendix 2: Computation of the Sectional Curvature for the Heisenberg Group}

We follow the plan:

\begin{itemize}
    \item[$(1)$] Compute the commutators $[X_i, X_j]$, $i = 1,2,3$, $j = 1,2,3$;

    \begin{align*}
    [X_1, X_2] = c_{12}^1X_1 + c_{12}^2X_2 + c_{12}^3X_3 = X_3;
    \end{align*}

    \begin{align*}
    [X_1, X_3] = c_{13}^1X_1 + c_{13}^2X_2 + c_{13}^3X_3 = 0;
    \end{align*}

    \begin{align*}
    [X_2, X_3] = c_{23}^1X_1 + c_{23}^2X_2 + c_{23}^3X_3 = 0.
    \end{align*}
    
    \item[$(2)$] Compute the Levi-Civita connection coefficients using the Koszul formula (\ref{koszul}), using the basis vector fields $X_1, X_2$, $X_3$;

    \begin{equation*}
    G_{ij}^k = 2g \left(D_{X_i}X_j,X_k \right) = -g\left( X_i,[X_j,X_k] \right) + g\left( X_j,[X_k,X_i] \right) + g\left( X_k,[X_i,X_j] \right).
    \end{equation*}

    \begin{align*}
    & D_{X_1}X_1 = \frac{1}{2}\left( 2c_{12}^1X_2 + 2c_{13}^1X_3 \right) = c_{12}^1X_2 + c_{13}^1X_3 = 0;\\
    & D_{X_1}X_2 =  \frac{1}{2}\left( -\left( - 2c_{12}^1 \right)X_1 +  \left( c_{23}^1 + c_{31}^2 + c_{12}^3 \right)X_3 \right) = c_{12}^1X_1 + \frac{1}{2}\left( c_{23}^1 + c_{31}^2 + c_{12}^3 \right)X_3 = \frac{1}{2}X_3;\\
    & D_{X_1}X_3 = \frac{1}{2}\left( -\left( -2c_{13}^1 \right)X_1 + \left( -c_{32}^1 + c_{21}^3 + c_{13}^2 \right)X_2 \right) = c_{13}^1X_1 - \frac{1}{2}\left( c_{32}^1 - c_{21}^3 - c_{13}^2 \right)X_2 = -\frac{1}{2}X_2; \\
    & D_{X_2}X_1 = \frac{1}{2}\left( \left( -c_{12}^2 - c_{22}^1 + c_{21}^2 \right)X_2 + \left( -c_{13}^2 - c_{32}^1 + c_{21}^3 \right)X_3 \right) = -c_{12}^2X_2 - \frac{1}{2}\left( c_{13}^2 + c_{32}^1 - c_{21}^3 \right)X_3 = -\frac{1}{2}X_3;\\
    & D_{X_2}X_2 = \frac{1}{2}\left( -\left( -c_{21}^2 + c_{12}^2 - c_{22}^1 \right)X_1 + \left( -c_{23}^2 + c_{32}^2 + c_{22}^3 \right)X_3 \right) = -c_{12}^2 X_1 - c_{23}^2 X_3 = 0;\\
    & D_{X_2}X_3 = \frac{1}{2}\left( -\left( -c_{31}^2 + c_{12}^3 - c_{23}^1 \right)X_1 + \left( -c_{32}^2 + c_{22}^3 + c_{23}^2 \right)X_2 \right) = \frac{1}{2}\left(c_{31}^2 - c_{12}^3 + c_{23}^1 \right) X_1 + c_{23}^2 X_2 = -\frac{1}{2}X_1;\\
    & D_{X_3}X_1 = \frac{1}{2}\left( \left( -c_{12}^3  - c_{23}^1 + c_{31}^2 \right)X_2 + \left( -c_{13}^3 - c_{33}^1 + c_{31}^3  \right)X_3 \right) = -\frac{1}{2}\left( c_{12}^3  + c_{23}^1 - c_{31}^2 \right)X_2 - c_{13}^3X_3 = -\frac{1}{2}X_2;\\
    & D_{X_3}X_2 = \frac{1}{2}\left( -\left( -c_{21}^3 + c_{13}^2 - c_{32}^1 \right)X_1 + \left( -c_{23}^3 + c_{33}^2 + c_{32}^3  \right)X_3 \right) = \frac{1}{2}\left( c_{21}^3-c_{13}^2 + c_{32}^1 \right)X_1 - c_{23}^3 X_3 = -\frac{1}{2}X_1;\\
    & D_{X_3}X_3 = \frac{1}{2}\left( -\left( -c_{31}^3 + c_{13}^3 - c_{33}^1 \right)X_1 + \left( -c_{32}^3 + c_{23}^3 + c_{33}^2 \right)X_2 \right) = -c_{13}^3X_1 + c_{23}^3 X_2 = 0.
    \end{align*}
    
    \item[$(3)$] Compute the Riemann curvature tensor $R_{v w}v$, where $v$, $w$ define some plane $P$ in the tangent space at point $q$;

    Expand the fields $v$, $w$ in the basis:
    \begin{equation*}
    v = v^1X_1 + v^2X_2 + v^3X_3, \quad w = w^1X_1 + w^2X_2 + w^3X_3.
    \end{equation*}

    Then

    \begin{align*}
    &  R_{vw}v = D_{[v,w]}v-[D_{v},D_{w}]v = R_{v,w,v}^1X_1 + R_{v,w,v}^2X_2 + R_{v,w,v}^3X_3. 
    \end{align*}

    Using the fact that $c_{ij}^k = 0$, except $c_{12}^3 = 1 = -c_{21}^3$, we obtain some simplification:

    \begin{align*}
    & R_{v,w,v}^k = \sum_{l,m=1}^3\left( w^lD_{v,v,l,m}^k - v^lD_{w,v,l,m}^k \right) + \sum_{l=1}^{3}B^l D_{v,l}^k = \sum_{l,m=1}^3\left( w^lD_{v,v,l,m}^k - v^lD_{w,v,l,m}^k \right) + \\
    & +\sum_{l=1}^3\left(\sum_{m=1}^2\sum_{m<j\leqslant 3}\left( v^mw^j -v^jw^m \right)c_{m,j}^l + \sum_{m=1}^3\left( \left( v^mX_m \right)w^l \right) - \sum_{m=1}^3w^m\left( X_m v^l \right) \right)D_{v,l}^k = \\
    & = \sum_{l,m=1}^3\left( w^lD_{v,v,l,m}^k - v^lD_{w,v,l,m}^k \right) + \left( v^1w^2 - v^2w^1 \right) D_{v,3}^k + \sum_{l, m=1}^3 \left( \left( \left( v^mX_m \right)w^l \right) - w^m\left( X_m v^l \right) \right) D_{v,l}^k, \quad k = 1,2,3.
    \end{align*}

    Next, we compute $D_{v,l}^k$, $D_{w,v,l,m}^k$, $D_{v,v,l,m}^k$, $l,\ m,\ k = 1,2,3$ in general form, and then for the case of constant coordinate functions $v^i$, $w^j$, $i,\ j = 1,2,3$:
     \begin{align*}
    & D_{v,1}^1 = \left( X_1 v^1 \right),\ D_{v,1}^2 = - \frac{v^3}{2} + \left( X_1 v^2 \right),\ D_{v,1}^3 = \frac{v^2}{2} + \left( X_1 v^3 \right);\\
    & D_{v,2}^1 = -\frac{v^3}{2} + \left( X_2 v^1 \right),\ D_{v,2}^2 = \left( X_2 v^2 \right),\ D_{v,2}^3 = -\frac{v^1}{2} + \left( X_2 v^3 \right),\\
    & D_{v,3}^1 = -\frac{v^2}{2} + \left( X_3 v^1 \right),\ D_{v,3}^2 = -\frac{v^1}{2} + \left( X_3 v^2 \right),\ D_{v,3}^3 = \left( X_3 v^3 \right);
    \end{align*}
    
    \begin{align*}
    & D_{w,v,1,1}^1 = \left( X_1 w^1 D_{v,1}^1 \right),\quad D_{w,v,1,1}^2 = - \frac{w^1D_{v,1}^3}{2} + \left( X_1 w^1 D_{v,1}^2  \right), \quad D_{w,v,1,1}^3 =  \frac{w^1D_{v,1}^2}{2} + \left( X_1 w^1 D_{v,1}^3 \right);\\
    & D_{w,v,2,1}^1 = - \frac{w^1D_{v,1}^3}{2} + \left( X_2 w^1 D_{v,1}^1 \right),\quad D_{w,v,2,1}^2 = \left( X_2 w^1 D_{v,1}^2  \right), \quad D_{w,v,2,1}^3 = -\frac{w^1D_{v,1}^1}{2} + \left( X_2 w^1 D_{v,1}^3 \right);\\
    & D_{w,v,3,1}^1 = -\frac{w^1D_{v,1}^2}{2} + \left( X_3 w^1 D_{v,1}^1 \right),\quad D_{w,v,3,1}^2 = -\frac{w^1D_{v,1}^1}{2} + \left( X_3 w^1 D_{v,1}^2  \right), \quad D_{w,v,3,1}^3 =  \left( X_3 w^1 D_{v,1}^3 \right);\\
    & D_{w,v,1,2}^1 = \left( X_1 w^2 D_{v,2}^1 \right),\quad D_{w,v,1,2}^2 = - \frac{w^2D_{v,2}^3}{2} + \left( X_1 w^2 D_{v,2}^2  \right), \quad  D_{w,v,1,2}^3 = \frac{w^2D_{v,2}^2}{2} + \left( X_1 w^2 D_{v,2}^3 \right);\\
    & D_{w,v,2,2}^1 = - \frac{w^2D_{v,2}^3}{2} + \left( X_2 w^2 D_{v,2}^1 \right),\quad D_{w,v,2,2}^2 = \left( X_2 w^2 D_{v,2}^2  \right), \quad D_{w,v,2,2}^3 = -\frac{w^2D_{v,2}^1}{2} + \left( X_2 w^2 D_{v,2}^3 \right);\\
    & D_{w,v,3,2}^1 = -\frac{w^2D_{v,2}^2}{2} + \left( X_3 w^2 D_{v,2}^1 \right),\quad D_{w,v,3,2}^2 = -\frac{w^2D_{v,2}^1}{2} + \left( X_3 w^2 D_{v,2}^2  \right), \quad D_{w,v,3,2}^3 = \left( X_3 w^2 D_{v,2}^3 \right);\\
    & D_{w,v,1,3}^1 = \left( X_1 w^3 D_{v,3}^1 \right),\quad D_{w,v,1,3}^2 = - \frac{w^3D_{v,3}^3}{2} + \left( X_1 w^3 D_{v,3}^2  \right), \quad D_{w,v,1,3}^3 =  \frac{w^3D_{v,3}^2}{2} + \left( X_1 w^3 D_{v,3}^3 \right);\\
    & D_{w,v,2,3}^1 = - \frac{w^3D_{v,3}^3}{2} + \left( X_2 w^3 D_{v,3}^1 \right),\quad D_{w,v,2,3}^2 = \left( X_2 w^3 D_{v,3}^2  \right), \quad D_{w,v,2,3}^3 = -\frac{w^3D_{v,3}^1}{2} + \left( X_2 w^3 D_{v,3}^3 \right); \\
    & D_{w,v,3,3}^1 = -\frac{w^3D_{v,3}^2}{2} + \left( X_3 w^3 D_{v,3}^1 \right),\quad D_{w,v,3,3}^2 = -\frac{w^3D_{v,3}^1}{2} + \left( X_3 w^3 D_{v,3}^2  \right), \quad D_{w,v,3,3}^3 = \left( X_3 w^3 D_{v,3}^3 \right).
    \end{align*}

    \begin{align*}
    & D_{v,v,1,1}^1 = \left( X_1 v^1 D_{v,1}^1 \right),\quad D_{v,v,1,1}^2 = - \frac{v^1D_{v,1}^3}{2} + \left( X_1 v^1 D_{v,1}^2  \right), \quad D_{v,v,1,1}^3 =  \frac{v^1D_{v,1}^2}{2} + \left( X_1 v^1 D_{v,1}^3 \right);\\
    & D_{v,v,2,1}^1 = - \frac{v^1D_{v,1}^3}{2} + \left( X_2 v^1 D_{v,1}^1 \right),\quad D_{v,v,2,1}^2 = \left( X_2 v^1 D_{v,1}^2  \right), \quad D_{v,v,2,1}^3 = -\frac{v^1D_{v,1}^1}{2} + \left( X_2 v^1 D_{v,1}^3 \right);\\
    & D_{v,v,3,1}^1 = -\frac{v^1D_{v,1}^2}{2} + \left( X_3 v^1 D_{v,1}^1 \right),\quad D_{v,v,3,1}^2 = -\frac{v^1D_{v,1}^1}{2} + \left( X_3 v^1 D_{v,1}^2  \right), \quad D_{v,v,3,1}^3 =  \left( X_3 v^1 D_{v,1}^3 \right);\\
    & D_{v,v,1,2}^1 = \left( X_1 v^2 D_{v,2}^1 \right),\quad D_{v,v,1,2}^2 = - \frac{v^2D_{v,2}^3}{2} + \left( X_1 v^2 D_{v,2}^2  \right), \quad  D_{v,v,1,2}^3 = \frac{v^2D_{v,2}^2}{2} + \left( X_1 v^2 D_{v,2}^3 \right);\\
    & D_{v,v,2,2}^1 = - \frac{v^2D_{v,2}^3}{2} + \left( X_2 v^2 D_{v,2}^1 \right),\quad D_{v,v,2,2}^2 = \left( X_2 v^2 D_{v,2}^2  \right), \quad D_{v,v,2,2}^3 = -\frac{v^2D_{v,2}^1}{2} + \left( X_2 v^2 D_{v,2}^3 \right);\\
    & D_{v,v,3,2}^1 = -\frac{v^2D_{v,2}^2}{2} + \left( X_3 v^2 D_{v,2}^1 \right),\quad D_{v,v,3,2}^2 = -\frac{v^2D_{v,2}^1}{2} + \left( X_3 v^2 D_{v,2}^2  \right), \quad D_{v,v,3,2}^3 = \left( X_3 v^2 D_{v,2}^3 \right);\\
    & D_{v,v,1,3}^1 = \left( X_1 v^3 D_{v,3}^1 \right),\quad D_{v,v,1,3}^2 = - \frac{v^3D_{v,3}^3}{2} + \left( X_1 v^3 D_{v,3}^2  \right), \quad D_{v,v,1,3}^3 =  \frac{v^3D_{v,3}^2}{2} + \left( X_1 v^3 D_{v,3}^3 \right);\\
    & D_{v,v,2,3}^1 = - \frac{v^3D_{v,3}^3}{2} + \left( X_2 v^3 D_{v,3}^1 \right),\quad D_{v,v,2,3}^2 = \left( X_2 v^3 D_{v,3}^2  \right), \quad D_{v,v,2,3}^3 = -\frac{v^3D_{v,3}^1}{2} + \left( X_2 v^3 D_{v,3}^3 \right); \\
    & D_{v,v,3,3}^1 = -\frac{v^3D_{v,3}^2}{2} + \left( X_3 v^3 D_{v,3}^1 \right),\quad D_{v,v,3,3}^2 = -\frac{v^3D_{v,3}^1}{2} + \left( X_3 v^3 D_{v,3}^2  \right), \quad D_{v,v,3,3}^3 = \left( X_3 v^3 D_{v,3}^3 \right).
    \end{align*}

    Next, we assume $v^i$, $w^i$ are constant. It follows that all derivatives $\left( X_i v^j \right)$, $\left( X_i w^j \right)$ are zero.

    We obtain

    \begin{align*}
    & D_{v,1}^1 = 0,\ D_{v,1}^2 = - \frac{v^3}{2},\ D_{v,1}^3 = \frac{v^2}{2};\\
    & D_{v,2}^1 = -\frac{v^3}{2},\ D_{v,2}^2 = 0,\ D_{v,2}^3 = -\frac{v^1}{2},\\
    & D_{v,3}^1 = -\frac{v^2}{2},\ D_{v,3}^2 = -\frac{v^1}{2},\ D_{v,3}^3 = 0;
    \end{align*}

    and also
    
    \begin{align*}
    & D_{w,v,1,1}^1 = 0,\quad D_{w,v,1,1}^2 = - \frac{w^1}{2}\frac{v^2}{2}, \quad 
    D_{w,v,1,1}^3 = - \frac{v^3}{2} \frac{w^1}{2};\\
    & D_{w,v,2,1}^1 = - \frac{w^1}{2}\frac{v^2}{2},\quad D_{w,v,2,1}^2 = 0, \quad D_{w,v,2,1}^3 = 0;\\
    & D_{w,v,3,1}^1 = \frac{v^3}{2}\frac{w^1}{2},\quad D_{w,v,3,1}^2 = 0, \quad 
    D_{w,v,3,1}^3 = 0;\\
    & D_{w,v,1,2}^1 = 0,\quad D_{w,v,1,2}^2 = \frac{v^1}{2} \frac{w^2}{2}, \quad  
    D_{w,v,1,2}^3 = 0;\\
    & D_{w,v,2,2}^1 = \frac{v^1}{2} \frac{w^2}{2},\quad D_{w,v,2,2}^2 = 0, \quad D_{w,v,2,2}^3 = \frac{w^2}{2}\frac{v^3}{2};\\
    & D_{w,v,3,2}^1 = 0,\quad D_{w,v,3,2}^2 = \frac{w^2}{2}\frac{v^3}{2}, \quad 
    D_{w,v,3,2}^3 = 0;\\
    & D_{w,v,1,3}^1 = 0,\quad D_{w,v,1,3}^2 = 0, \quad 
    D_{w,v,1,3}^3 = -\frac{v^1}{2} \frac{w^3}{2};\\
    & D_{w,v,2,3}^1 = 0,\quad D_{w,v,2,3}^2 = 0, \quad D_{w,v,2,3}^3 = \frac{v^2}{2}\frac{w^3}{2}; \\
    & D_{w,v,3,3}^1 = \frac{v^1}{2}\frac{w^3}{2},\quad D_{w,v,3,3}^2 = \frac{v^2}{2}\frac{w^3}{2}, \quad 
    D_{w,v,3,3}^3 = 0.
    \end{align*}

    \begin{align*}
    & D_{v,v,1,1}^1 = 0,\quad D_{v,v,1,1}^2 = - \frac{v^1}{2}\frac{v^2}{2}, \quad 
    D_{v,v,1,1}^3 =  - \frac{v^3}{2}\frac{v^1}{2};\\
    & D_{v,v,2,1}^1 = - \frac{v^1}{2}\frac{v^2}{2},\quad D_{v,v,2,1}^2 = 0, \quad D_{v,v,2,1}^3 = 0;\\
    & D_{v,v,3,1}^1 = \frac{v^3}{2}\frac{v^1}{2},\quad D_{v,v,3,1}^2 = 0, \quad 
    D_{v,v,3,1}^3 =  0;\\
    & D_{v,v,1,2}^1 = 0,\quad D_{v,v,1,2}^2 = \frac{v^1}{2} \frac{v^2}{2}, \quad  
    D_{v,v,1,2}^3 = 0;\\
    & D_{v,v,2,2}^1 = \frac{v^1}{2}\frac{v^2}{2},\quad D_{v,v,2,2}^2 = 0, \quad D_{v,v,2,2}^3 = \frac{v^2}{2}\frac{v^3}{2};\\
    & D_{v,v,3,2}^1 = 0,\quad D_{v,v,3,2}^2 = \frac{v^2}{2}\frac{v^3}{2}, \quad 
    D_{v,v,3,2}^3 = 0;\\
    & D_{v,v,1,3}^1 = 0,\quad D_{v,v,1,3}^2 = 0, \quad 
    D_{v,v,1,3}^3 = -\frac{v^1}{2} \frac{v^3}{2};\\
    & D_{v,v,2,3}^1 = 0,\quad D_{v,v,2,3}^2 = 0, \quad D_{v,v,2,3}^3 = \frac{v^2}{2}\frac{v^3}{2}; \\
    & D_{v,v,3,3}^1 = \frac{v^1}{2}\frac{v^3}{2},\quad D_{v,v,3,3}^2 = \frac{v^2}{2}\frac{v^3}{2}, \quad D_{v,v,3,3}^3 = 0.
    \end{align*}
    
    \item[$(4)$] Compute $Q(v,w)$, where $v$, $w$ define some plane $P$ in the tangent space at point $q$;

    \begin{equation*}
    Q(v,w) = \left( -\left( v^1 \right)^2 +  \sum_{k=2}^3\left( v^k \right)^2 \right) \cdot \left( -\left( w^1 \right)^2 +  \sum_{l=2}^3\left( w^l \right)^2 \right) - \left( -v^1w^1 + \sum_{l=2}^3v^lw^l\right)^2.
    \end{equation*}
    
    \item[$(5)$] Obtain the sectional curvature $K$ using the performed computations and formula (\ref{sectionalcurvature}).

    In the general case:

    \begin{equation*}
    K = \frac{g\left( R_{v,w}v, w \right)}{Q(v,w)} = \frac{-w^1R_{v,w,v}^1 + w^2R_{v,w,v}^2 + w^3R_{v,w,v}^3}{\left( -\left( v^1 \right)^2 +  \sum_{k=2}^3\left( v^k \right)^2 \right) \cdot \left( -\left( w^1 \right)^2 +  \sum_{l=2}^3\left( w^l \right)^2 \right) - \left( -v^1w^1 + \sum_{l=2}^3v^lw^l\right)^2}.
    \end{equation*}

    where

    \begin{equation*}
   R_{v,w,v}^k = \sum_{l,m=1}^3\left( w^lD_{v,v,l,m}^k - v^lD_{w,v,l,m}^k \right) + \sum_{l, m=1}^3 \left( \left( \left( v^mX_m \right)w^l \right) - w^m\left( X_m v^l \right) \right) D_{v,l}^k + \left( v^1w^2 - v^2w^1 \right) D_{v,3}^k, \quad k = 1,2,3.
    \end{equation*}

    For the case of constant coordinate functions, let us compute $R_{v,w,v}^k$, $k=1,2,3$:

    \begin{align*}
    & R_{v,w,v}^1 = \sum_{l,m=1}^3\left( w^lD_{v,v,l,m}^1 - v^lD_{w,v,l,m}^1 \right) + \left( v^1w^2 - v^2w^1 \right) D_{v,3}^1 =  \sum_{m=1}^3\left( w^1D_{v,v,1,m}^1 - v^1D_{w,v,1,m}^1 \right) + \\
    & + \sum_{m=1}^3\left( w^2D_{v,v,2,m}^1 - v^2D_{w,v,2,m}^1 \right) + \sum_{m=1}^3\left( w^3D_{v,v,3,m}^1 - v^3D_{w,v,3,m}^1 \right) - \frac{v^2}{2}\left( v^1w^2 - v^2w^1 \right) = \\
    & = w^2\left( - \frac{v^1}{2}\frac{v^2}{2} + \frac{v^1}{2}\frac{v^2}{2} \right) + w^3\left( \frac{v^3}{2}\frac{v^1}{2} + \frac{v^1}{2}\frac{v^3}{2} \right) - v^2\left( - \frac{w^1}{2}\frac{v^2}{2} + \frac{v^1}{2} \frac{w^2}{2} \right) - v^3\left( \frac{v^3}{2}\frac{w^1}{2} + \frac{v^1}{2}\frac{w^3}{2} \right) - \frac{v^2}{2}\left( v^1w^2 - v^2w^1 \right) = \\
    & = \frac{w^1}{4}\left( 3v^2v^2 - v^3v^3 \right) - w^2\frac{3v^1v^2}{4} + w^3\frac{v^1v^3}{4}; \\
    & R_{v,w,v}^2 = \sum_{l,m=1}^3\left( w^lD_{v,v,l,m}^2 - v^lD_{w,v,l,m}^2 \right) + \left( v^1w^2 - v^2w^1 \right) D_{v,3}^2 = \sum_{m=1}^3\left( w^1D_{v,v,1,m}^2 - v^1D_{w,v,1,m}^2 \right) + \\
    & + \sum_{m=1}^3\left( w^2D_{v,v,2,m}^2 - v^2D_{w,v,2,m}^2 \right) + \sum_{m=1}^3\left( w^3D_{v,v,3,m}^2 - v^3D_{w,v,3,m}^2 \right) -\frac{v^1}{2}\left( v^1w^2 - v^2w^1 \right)  = \\
    & = w^1\left( - \frac{v^1}{2}\frac{v^2}{2} + \frac{v^1}{2} \frac{v^2}{2} \right) +  w^3\left( \frac{v^2}{2}\frac{v^3}{2} + \frac{v^2}{2}\frac{v^3}{2} \right) - v^1\left( - \frac{w^1}{2}\frac{v^2}{2} + \frac{v^1}{2} \frac{w^2}{2} \right) - v^3\left(\frac{w^2}{2}\frac{v^3}{2} + \frac{v^2}{2}\frac{w^3}{2} \right) -\frac{v^1}{2}\left( v^1w^2 - v^2w^1 \right)  =  \\
    & = w^1\frac{3v^1v^2}{4} - \frac{w^2}{4}\left( 3v^1v^1 + v^3v^3 \right) + w^3\frac{v^2v^3}{4}; \\
    & R_{v,w,v}^3 = \sum_{l,m=1}^3\left( w^lD_{v,v,l,m}^3 - v^lD_{w,v,l,m}^3 \right) + \left( v^1w^2 - v^2w^1 \right) D_{v,3}^3 = \sum_{m=1}^3\left( w^1D_{v,v,1,m}^3 - v^1D_{w,v,1,m}^3 \right) + \\
    & + \sum_{m=1}^3\left( w^2D_{v,v,2,m}^3 - v^2D_{w,v,2,m}^3 \right) + \sum_{m=1}^3\left( w^3D_{v,v,3,m}^3 - v^3D_{w,v,3,m}^3 \right) = \\
    & = w^1\left( - \frac{v^3}{2}\frac{v^1}{2} -\frac{v^1}{2} \frac{v^3}{2} \right) + w^2\left( \frac{v^2}{2}\frac{v^3}{2} + \frac{v^2}{2}\frac{v^3}{2} \right) - v^1\left( - \frac{v^3}{2} \frac{w^1}{2} -\frac{v^1}{2} \frac{w^3}{2} \right) - v^2\left( \frac{w^2}{2}\frac{v^3}{2} + \frac{v^2}{2}\frac{w^3}{2} \right) = \\
    & = -w^1\frac{v^1v^3}{4} + w^2\frac{v^2v^3}{4} + \frac{w^3}{4}\left( v^1v^1 - v^2v^2 \right),
    \end{align*}

    and the numerator in the curvature formula:

    \begin{align*}
    & -w^1R_{v,w,v}^1 + w^2R_{v,w,v}^2 + w^3R_{v,w,v}^3 = -w^1\left( \frac{w^1}{4}\left( 3v^2v^2 - v^3v^3 \right) - w^2\frac{3v^1v^2}{4} + w^3\frac{v^1v^3}{4} \right) + \\
    & + w^2\left( w^1\frac{3v^1v^2}{4} - \frac{w^2}{4}\left( 3v^1v^1 + v^3v^3 \right) + w^3\frac{v^2v^3}{4} \right) + w^3\left( -w^1\frac{v^1v^3}{4} + w^2\frac{v^2v^3}{4} + \frac{w^3}{4}\left( v^1v^1 - v^2v^2 \right) \right) = \\
    & = -\frac{w^1w^1}{4}\left( 3v^2v^2 - v^3v^3 \right) -\frac{w^2w^2}{4}\left( 3v^1v^1 + v^3v^3 \right) + \frac{w^3w^3}{4}\left( v^1v^1 - v^2v^2 \right) + w^1w^2\frac{3v^1v^2}{2} - w^1w^3\frac{v^1v^3}{2} + w^2w^3\frac{v^2v^3}{2}.
    \end{align*}

    Then the final formula looks like this:

    \begin{equation*}
    K = \frac{-\frac{w^1w^1}{4}\left( 3v^2v^2 - v^3v^3 \right) -\frac{w^2w^2}{4}\left( 3v^1v^1 + v^3v^3 \right) + \frac{w^3w^3}{4}\left( v^1v^1 - v^2v^2 \right) + w^1w^2\frac{3v^1v^2}{2} - w^1w^3\frac{v^1v^3}{2} + w^2w^3\frac{v^2v^3}{2}}{\left( -\left( v^1 \right)^2 +  \sum_{k=2}^3\left( v^k \right)^2 \right) \cdot \left( -\left( w^1 \right)^2 +  \sum_{l=2}^3\left( w^l \right)^2 \right) - \left( -v^1w^1 + v^2w^2 + v^3w^3 \right)^2}.
    \end{equation*}
    
\end{itemize}

\end{document}